\documentclass[11pt]{amsart}
\usepackage{epic,eepic}
\usepackage{amsmath}
\usepackage{amssymb}
\input psfig  

\catcode`@=11 \@mparswitchfalse  

\date{July, 1998} 

\newcommand{\bi}{\bf\itshape}  

\numberwithin{equation}{section}   

\newcommand{\la}{\langle}	   	
\newcommand{\ra}{\rangle}          	

\renewcommand{\(}{\left(}
\renewcommand{\)}{\right)}

\newcommand{\R}{{\mathbf R}}		 
\newcommand{\e}{\varepsilon}             
\newcommand{\nothing}{\varnothing}       

\newcommand{\note}[1]{}                    

\newcommand{\mnote}[1]{}

\newcommand{\ralph}{{\bf *RH*} }  

\theoremstyle{plain}

\newtheorem{theorem}{Theorem}
\newtheorem{defn}{Definition}

\swapnumbers 

\newtheorem{thm}{Theorem}[section]
\newtheorem{lemma}[thm]{Lemma}
\newtheorem{prop}[thm]{Proposition}
\newtheorem{cor}[thm]{Corollary}

\theoremstyle{definition}

\theoremstyle{remark}
\newtheorem{remark}[thm]{Remark}

\newtheorem{example}[thm]{Example}

\renewcommand{\phi}{\varphi}				
\newcommand{\cd}{,\dots,}	
\newcommand{\ol}{\overline}
\newcommand{\dist}{\operatorname{dist}\nolimits}
\newcommand{\cn}{\colon}	

\newcommand{\conv}{\operatorname{Co}\nolimits}	
\newcommand{\X}{{\mathbf X}}		
\newcommand{\Nat}{{\mathbf N}}			
\newcommand{\supp}{\operatorname{supp}\nolimits}
\newcommand{\Q}{{\mathbf Q}}			
\newcommand{\dy}{{\mathcal D}}			
\newcommand{\simp}{{\mathcal S}}	
\newcommand{\maxside}{{\mathcal M}}	

\newcommand{\fp}[1]{\boldsymbol{\{}#1\boldsymbol{\}} }	
	
\title[Approximately Convex Functions and the Size of Convex
Hulls]{Extremal Approximately Convex Functions and 
Estimating the Size of Convex Hulls}

\author{S. J. Dilworth}
\address{\null\hskip-1\parindent Department of Mathematics\newline
University of South Carolina\newline
Columbia, S.C. 29208, USA\newline 
\null\newline
{\tt dilworth\char'100math.sc.edu}\newline 
{\tt howard\char'100math.sc.edu}\newline
{\tt roberts\char'100math.sc.edu}
}

\author{Ralph Howard}

\author{James W. Roberts}

\thanks{The research of the second author was supported in part from ONR Grant
N00014-90-J-1343 and ARPA-DEPSCoR Grant DAA04-96-1-0326}

\begin{document}

\begin{abstract} 
A real valued function $f$ defined on a convex $K$ is an {\bi approximately
convex function\/} iff it satisfies
$$
f\Bigl(\frac{x+y}{2}\Bigr)\le \frac{f(x)+f(y)}{2}+1.
$$
A thorough study of approximately convex functions is made.  The principal
results are a sharp universal upper bound for lower semi-continuous
approximately convex functions that vanish on the vertices of a simplex and
an explicit description of the unique largest bounded
approximately convex function~$E$ vanishing on the vertices of a simplex.

A set $A$ in a normed space is an {\bi approximately convex set\/} iff
for all $a,b\in A$ the distance of the midpoint $(a+b)/2$ to $A$ is
$\le 1$.  The bounds on approximately convex functions are used to show
that in $\R^n$ with the Euclidean norm, for any approximately convex set $A$,
any point $z$ of the convex hull of $A$ is at a distance of at most
$[\log_2(n-1)]+1+(n-1)/2^{[\log_2(n-1)]}$ from $A$.  Examples are
given to show this is the sharp bound.  Bounds for general norms on
$\R^n$ are also given.
\end{abstract}

\maketitle

\tableofcontents

\section{Introduction}

The problem motivating this paper is the following: given a set $A$ in
$\R^n$, estimate the size of the convex hull $\conv(A)$ of $A$ in terms
of geometric properties of $A$.  To do this we assume that $\R^n$ is
equipped with a norm $\|\cdot\|$.  Then a first step in constructing
the convex hull of $A$ is to add all the midpoints of segments joining
points of $A$.  The size of $\conv(A)$ can be estimated in terms of
this first step.

Our main result gives the sharp constants in this estimate for
$n$-dimensional Euclidean spaces and it provides an estimate for the
constants for general $n$-dimensional normed spaces which is accurate
to within $2/n$.  (Here $\dist(x,A):=\inf\{\|x-a\|:a\in A\}$ is the
distance of the point $x$ from the set $A$.)

\begin{theorem}\label{midpoint}  
If $(\R^n,\|\cdot\|)$ is an $n$-dimensional normed linear space then
there is a constant $C_{\|\cdot\|}$, depending on the norm $\|\cdot\|$,
so that if $A\subset \R^n$ satisfies
\begin{equation}\label{mid-dist-delta}
a_0,a_1\in A\quad \text{implies}\quad 
	\dist\Big(\frac{a_0+a_1}{2},A\Big)\le \delta,
\end{equation}
then
$$
z\in \conv(A)\quad \text{implies} \quad \dist(z,A)\le C_{\|\cdot\|}\delta.
$$
Letting $[\,\cdot\,]$ be the greatest integer function, the sharp
constant $C_{\|\cdot\|}$ satisfies
$$
[\log_2(n-1)]+1+(n-1)/2^{[\log_2(n-1)]}\le C_{\|\cdot\|}
	\le [\log_2(n)]+1+n/2^{[\log_2(n)]}
$$
(this holds for all norms) and the sharp constant when $\|\cdot\|$ is
the Euclidean norm is $C_{\|\cdot\|}=[\log_2(n-1)]+1+(n-1)/2^{[\log_2(n-1)]}$.
\end{theorem}

The upper bound $C_{\|\cdot\|}\le 2\lceil \log_2(n+1) \rceil$ (where
$\lceil\cdot \rceil$ is the ceiling function) is implicit in the
paper~\cite[Props~3.3 and 3.4]{Casini-Papini:almost-convex} of Casini
and Papini.

For bounded sets this can be given a concise restatement in terms of
the Hausdorff distance between sets.  Recall that if $A,B\subset \R^n$
are bounded then the {\bi Hausdorff distance\/}, $d_H(A,B)$,
between $A$ and $B$ is the infimum of the numbers $r$ so that every
point of $A$ is within a distance $r$ of a point of $B$ and every
point of $B$ is within distance $r$ of a point of $A$.
Define numbers
$\kappa(n)$ for $n\ge 0$ by $\kappa(0):=0$ and
$$
\kappa(n):= [\log_2(n)]+1+n/2^{[\log_2(n)]}
$$
for $n\ge 1$.  The collection of midpoints of segments joining pairs of
points of $A$ is $\frac12(A+A)=\{(a+b)/2:a,b\in A\}$.  Then
Theorem~\ref{midpoint} can be restated as
$$
d_H(\conv(A),A)\le C_{\|\cdot\|}d_H\Big(\frac12(A+A),A\Big)
$$
where the sharp constant $C_{\|\cdot\|}$ satisfies
$$
\kappa(n-1)\le C_{\|\cdot\|}\le \kappa(n)
$$
and $C_{\|\cdot\|}=\kappa(n-1)$ when $\|\cdot\|$ is the Euclidean
norm.  (Allowing $+\infty$ for a value of $d_H(A,B)$ this also holds
for unbounded sets).
If $A$ is a finite set with $N$ points, then
$d_H(A,\frac12(A+A))=\max_{b,c\in A}\min_{a\in A}\|a-\frac12(b+c)\|$
which can be computed in $O(N^3)$ operations.  Thus for finite sets
Theorem~\ref{midpoint} allows  estimation of $d_H(A,\conv(A))$ in
polynomial time.

For general norms obtaining the lower bound $\kappa(n-1)\le
C_{\|\cdot\|}$ is more difficult than the upper bound and
involves construction of some interesting geometric objects, the
extremal approximately convex functions.  To describe these we first
make a couple of definitions.  The following is motivated by taking
$\delta=1$ in the hypothesis of Theorem~\ref{midpoint}.

\begin{defn}\label{al-convex-set}
Let $(\X,\|\cdot\|)$ be a normed space.  Then a subset $A\subset \X$
is an {\bi approximately convex set\/} iff for all $a,b\in A$
$$
\dist\Bigl(\frac12(a+b),A\Bigr)\le 1.
$$
\vskip-15pt \qed
\end{defn}
If $A$ is an approximately convex set then the function $h(x)=\dist(x,A)$
(the distance of $x$ from $A$) will satisfy a weak form of the
inequality satisfied by a convex function.  We isolate this 
property:
\begin{defn}\label{al-convex-fcn}
Let $E$ be a convex set in the normed space $(\X,\|\cdot\|)$.  Then a
function $h\cn E\to \R$ is an {\bi approximately convex function\/}
iff for all $a,b\in E$
$$
h\Bigl(\frac{a+b}2\Bigr)\le \frac{h(a)+h(b)}{2}+1.
$$
\vskip-15pt \qed
\end{defn}
\noindent
(Strictly speaking this should be ``approximately midpoint convex'' or
``approximately Jensen convex'' but for the sake of brevity we will
use ``approximately convex''.)   Let $\Delta_n:=\{(\alpha_0\cd
\alpha_n): \alpha_k\ge 0,\ \sum_{k=0}^n\alpha_k\}$ be the standard
$n$-dimensional simplex.  Then the result leading to the lower bounds
on $C_{\|\cdot\|}$ is the explicit computation of the extremal
approximately convex function on the simplex.  

\begin{theorem}\label{E-intro}
There is an approximately convex $E\cn \Delta_n\to \R$ which vanishes
on the vertices of $\Delta_n$ with the following properties:
\begin{enumerate}
	\item \label{E-intro-char}
	If $h$ is a bounded (or Borel-measurable) approximately
	convex function on $\Delta_n$ which takes non-positive values 
	on the vertices, then $h(x) \le E(x)$ for $x\in \Delta_n$.
	\item $E$ achieves its maximum value of $\kappa(n)$.
	\item $E$ is lower semi-continuous.
\end{enumerate}
The property~\ref{E-intro-char} characterizes $E$ uniquely.  Moreover
$E$ is given concretely in terms of an elementary infinite sum (see
equations (\ref{H-frac-series}) and (\ref{E-defn})).
\end{theorem}

The examples showing the lower bounds on $C_{\|\cdot\|}$ in
Theorem~\ref{midpoint} are sharp and are constructed from the graph of
$E$.  The lower semi-continuity of $E$ and the fact that $E$ has $\kappa(n)$ as
its maximum are important in these constructions.  We note the mere
existence of $E$ (which follows from abstract considerations) is less
important than the fact that $E$ is given explicitly in a relatively
simple form (cf.~\S\ref{sec:E} and Figure~\ref{E2d}).
\medskip

We now give a more detailed description of our
results. In~\S\ref{sec:bounds} we give upper bounds on approximately
convex functions which are locally bounded from above.  Motivated by
Perron's method in the theory of harmonic functions 
in~\S\ref{sec:lsc} we show that given a compact convex set
$K\subset \R^n$ with extreme points $V$ and a uniformly continuous
function $\phi\cn V\to \R$ then there is a a unique extremal bounded
approximately convex function $E_{K,\phi}$ on $K$ which agrees with
$\phi$ on $V$; moreover, $E_{K,\phi}$ is realized (as in Perron's
method) as the pointwise supremum of all bounded approximately convex
functions on $K$ which agree with $\phi$ on $V$.  The function
$E_{K,\phi}$ is lower semi-continuous, characterized by a mean-value
property, and satisfies a certain maximum principle.

~\S\ref{sec:H} and~\S\ref{sec:E} contain a description of the
extremal approximately convex function~$E$ on the simplex and proofs
of the properties of~$E$ listed in Theorem~\ref{E-intro}.  In
\S~\ref{sec:polytopes}, we determine the extremal function
$E_{K,\phi}$ when $K$ is a convex polytope.  A stability theorem with
sharp constants for approximately convex functions of the type first
given by Hyers and Ulam~\cite{Hyers-Ulam:convex} is given
in~\S\ref{sec:Hyers-Ulam}.  This states that an approximately convex
function can be approximated in the uniform norm by a convex function
with error only depending on the dimension of the domain.  The example
showing the constants are sharp is the extremal function~$E$.  The
rest of Section~\ref{sec:functions} gives various other properties and
examples of approximately convex functions.

Section~\ref{sec:sets} gives the proof of Theorem~\ref{midpoint} and some
of its extensions and refinements.  The first two sections give the
upper and lower bounds $\kappa(n-1)\le C_{\|\cdot\|}\le \kappa(n)$ for
general norms.  The upper bound follows from the  general upper
bounds on approximately convex functions and the lower bound uses
properties of the extremal approximately convex function $E$ on
$\Delta_n$.
The proof that $C_{\|\cdot\|} = \kappa(n-1)$ in the
Euclidean case is given in \S\ref{sec:euclidean}.  This requires some
(hopefully interesting) geometrical arguments in addition
to Theorem~\ref{E-intro}.  Finally, we prove that $C_{\|\cdot\|}=2$
for all two-dimensional norms.  This argument is somewhat ad hoc and
does not appear to extend to higher dimensions.

\section{Approximately Convex Functions}\label{sec:functions}

We first relate approximately convex functions to approximately convex
sets.
\begin{prop} \label{set=fcn}
Let $(\X,\|\cdot\|)$ be a normed space, $A\subset \X$, and define
$h(x):=\dist(x,A)$.  Then $A$ is an approximately convex set if and only if
$h$ is an approximately convex function.
\end{prop}

\begin{proof}
If $h(x)=\dist(x,A)$ is an approximately convex function it is clear
that $A$ is an approximately convex set.  Conversely if $A$ is an
approximately convex set, let $x_0,x_1\in \X$ and $\e>0$.  Choose
$a_0, a_1\in A$ so that $h(x_0)=\dist(x_0,A)\le \|x_0-a_0\|+\e$ and
$h(x_1)\le
\|x_1-a_1\|+\e$.  As $A$ is approximately convex
$\dist((a_0+a_1)/2,A)\le 1$.  Thus
\begin{align*}
\dist\Bigl(\frac{x_0+x_1}{2},A\Bigr)&\le
\Bigl\|\frac{x_0-a_0}{2}\Bigr\|+\Bigl\|\frac{x_1-a_1}{2}\Bigr\|+
\dist(\frac{a_0+a_1}{2},A)\\
&\le \frac{h(x_0)+h(x_1)}{2}+\e+1
\end{align*}
As $\e$ as arbitrary this completes the proof.  (This proof is implicit
in the paper of Casini and
Papini~\cite[Prop.~3.4]{Casini-Papini:almost-convex}.)
\end{proof}

\subsection{Bounds on approximately convex functions}\label{sec:bounds}  

The first bound  is an extension to approximately convex functions of a
standard result about  convex functions.

\begin{prop}\label{lower-bd}
Let $U\subseteq \R^n$ be a convex  set and $h\cn U\to \R$ be approximately
convex and bounded from above by $C$.  Then for any $x_0\in U$ and
$x\in U\cap (2x_0-U)$ (if $x_0$ is in the interior of $U$ this is a
neighborhood of $x_0$ in $U$) the inequality
$$
h(x)\ge 2h(x_0)-C-2
$$
holds, and so $h$ is bounded from below in $U\cap (2x_0-U)$.  Thus $h$
is bounded from below on compact subsets of the interior of $U$.
\end{prop}

\begin{proof}  Let $y=2x_0-x$.  Then $y\in U$ as $x\in (2x_0-U)$.
Also $x_0=(x+y)/2$.  Thus
$$
h(x_0)=h\Bigl(\frac{x+y}{2}\Bigr)\le 
	\frac{h(x)+h(y)}{2}+1\le\frac{h(x)+C}{2}+1.
$$
Solving this for $h(x)$ completes the proof.
\end{proof}

The following theorem is one of our main results.

\begin{thm}\label{thm:convex-bd}
Let $A\subset \R^n$ with convex hull $E=\conv(A)$.  Let $h\cn E\to \R$
be an approximately convex function which is bounded above and which
satisfies $h \le 0$ on $A$.  Then
$$
\sup_{x\in E} h(x)\le  [\log_2n]+1+\dfrac{n}{2^{[\log_2n]}}.
$$
Moreover this is the sharp upper bound (the sharpness follows from
Theorem~\ref{II3}). 
\end{thm}

\begin{remark} The assumption that $h$ is bounded above can not be
dropped.  For the relevant example see Example~\ref{ex:Q-linear} in
\S\ref{sec:examples} below.\qed
\end{remark}

Before giving the proof we give a name to the bounds in the Theorem
and show that they satisfy a recursion which is a main ingredient of
the proof.  Let $\kappa(0)=0$ and for $n\ge 1$ 
\begin{equation}\label{kappa-def}
\kappa(n)= [\log_2n]+1+\dfrac{n}{2^{[\log_2n]}}.
\end{equation}
This notation will be use throughout the rest of the paper.

\begin{prop}\label{prop:kappa-recur}
The sequence $\la \kappa(n)\ra_{k=0}^\infty$  satisfies the recursion
\begin{equation}\label{kappa-recur}
\kappa(n)=\max_{\begin{subarray}{c} n_1+n_2=n\\ \vphantom{I} n_1,n_2\ge
0\end{subarray}} \frac{\kappa(n_1)+\kappa(n_2)}{2} +1
\end{equation}
for $n\ge 1$.
\end{prop}

\mnote{\ralph Here and below replaced ``decreasing'' with ``monotone
decreasing''} 
\begin{lemma}\label{concave-seq}
Let $\la \alpha(i)\ra_{i=0}^m$ be a finite sequence on $\{0,1\cd m\}$
so that $\la \alpha(j)-\alpha(j-1)\ra_{j=1}^m$ is monotone decreasing
(that is the sequence is concave).  Then
$$
\max_{i+j=n}\frac{\alpha(i)+\alpha(j)}{2}+1=
\begin{cases} \alpha(n)+1=\dfrac{\alpha(n)+\alpha(n)}{2}+1,&m=2n;\\
	\dfrac{\alpha(n)+\alpha(n+1)}{2}+1,&m=2n+1.\end{cases}
$$
\end{lemma}

\begin{proof}
Let $\beta(i)=(\alpha(i)+\alpha(n-i))/{2}+1$.  Then the concavity of
$\la \alpha(i)\ra$ implies the sequence $\la \beta(i)\ra$ is also
concave.  Also $\beta(i)=\beta(n-i)$ so $\la \beta(i)\ra$ is
symmetric.  But a symmetric concave function takes on its maximum at
the center of its interval of definition.  Thus if $m=2n$ is even the
maximum is $\beta(n)=\alpha(n)+1$ and if $m=2n+1$ the maximum is
$\beta(n)=\beta(n+1) =(\alpha(n)+\alpha(n+1))/2+1$.
\end{proof}

\begin{proof}[Proof of Proposition~\ref{prop:kappa-recur}]
A calculation shows
$$
\kappa(2n)=\kappa(n)+1,\quad \kappa(2n+1)=\frac{\kappa(n)+\kappa(n+1)}{2}+1.
$$
(The second of these is most easily seen by writing $n=2^m+r$ where
$0\le r\le 2^m-1$.) But the sequence $\la
\kappa(n)-\kappa(n-1)\ra_{k=1}^\infty$ is monotone decreasing
so that an application of the last lemma completes the proof.
\end{proof}

\begin{proof}[Proof of Theorem~\ref{thm:convex-bd}]
Recalling the definition of $\kappa(n)$ we wish to show that
$\sup_{x\in E}h(x)\le \kappa(n)$ 
We use induction on~$n$ based on the recursion~(\ref{kappa-recur})
satisfied by $\kappa$.  The base case of $n=0$ is clear.  Suppose
$n\ge 1$ and assume that the assertion holds for all integers less
than~$n$.  If $x\in E$ then by Carath\'eodory's Theorem
(cf.~\cite[p.~3]{Schneider:convex}) there are $x_0\cd x_n\in A$ so
that $x\in \conv\{x_0\cd x_n\}$.  Thus without loss of generality we may
assume that $E=\conv\{x_0\cd x_n\}$.  Let
$M:=\sup_{x\in E} h(x)<\infty $
and let $\e>0$.  Suppose that $x=\sum_{k=0}^n\alpha_kx_k\in E$ (with
$\sum_{k=0}^n\alpha_k=1$ and $\alpha_k\ge 0$) and $h(x)\ge M-\e$.
By reordering the terms if necessary we may assume $\alpha_0\le
\alpha_1\le \dots\le\alpha_n$.  Note that 
$\alpha_0\le 1/(n+1)\le 1/2$.  Let
$n_1$ be the least integer so that
$$
\sum_{k=0}^{n_1}\alpha_k  > \frac12.
$$
Then $\sum_{k=0}^{n_1-1}\alpha_k\le \frac12$.
Set
$$
s=\frac12-\sum_{k=0}^{n_1-1}\alpha_k, \quad t=\alpha_{n_1}-s, 
$$
and let
$$
y=2\biggl(\sum_{k=0}^{n_1-1}x_k+sx_{n_1}\biggr), \quad z=2\biggl(tx_{n_1}+\sum_{n_1+1}^n\alpha_kx_k\biggr).
$$
Then $y\in E$ as $\sum_{k=0}^{n_1-1}\alpha_k+s=\frac12$. Likewise
$z\in E$.  In particular $y\in \conv\{x_0\cd x_{n_1}\}=:\Delta_1$ and
$z\in \conv\{ x_{n_1}\cd x_n\}=:\Delta_2$.  Then $\dim \Delta_1=n_1$
and $\dim \Delta_2=n-n_1=:n_2$.  Since $\alpha_0\le 1/2$ we have
$n_1\ge 1$.

If $n_1<n$ then $n_1,n_2<n$ and therefore by the induction hypothesis
and $x=\frac12(y+z)$, we have
\begin{align*}
M -\e& \le h(x)\le \frac{h(y)+h(z)}{2}+1
 \le \frac{\kappa(n_1)+\kappa(n_2)}{2}+1\\
&\le \kappa(n_1+n_2)=\kappa(n).
\end{align*}
Therefore $M\le \kappa(n)+\e$.  
This leaves the case  $n_1=n$.  Then $z=x_n\in A$ and thus
$h(z)=0$. Whence
$$
M-\e\le h(x)\le \frac{h(y)+h(z)}2+1=\frac{h(y)}{2}+1\le \frac{M}2+1.
$$
Solve this inequality for $M$ and use $2\le \kappa(n)$ to get $M\le
2(1+\e)\le \kappa(n)(1+\e)$. Combining the inequalities from the two
cases and letting $\e \searrow 0$ implies $M\le \kappa(n) $ and
completes the proof.
\end{proof}

\begin{remark}\label{kappa-size}
As many of our results will involve $\kappa(n)$ it is worth giving
some sharp bounds on $\kappa(n)$.  To do this extend $\kappa$ to the
positive reals by defining $\kappa(x)=[\log_2x]+1+x/{2^{[\log_2x]}}$.
Then for any integer $m$ we have $\kappa(2^m)=m+2=\log_2(2^m)+2$.  On
closed intervals $[2^m,2^{m+1}]$ the function $\kappa(x)$ is linear.
Thus $\kappa(x)$ is the continuous piecewise linear function on
$(0,\infty)$ with knots at $x=2^m$ and with $\kappa(x)=2+\log_2(x)$ at
the knots.  As the function $2+\log_2(x)$ is concave this implies
$\kappa(x) \le 2+\log_2(x)$.  On each of the intervals it is a
straightforward calculus exercise to find the maximum of
$(2+\log_2(x))-\kappa(x)$ on the interval $[2^m,2^{m+1}]$.  The result
is $(\ln(2)-\ln(\ln(2))-1)/\ln(2)\approx .08607133206$ (surprisingly
this is independent of which interval $[2^m,2^{m+1}]$ we are working
on).  This leads to the bounds
$$
1.913928+\log_2(n)<\kappa(n)\le 2+\log_2(n).
$$
\vskip-15pt\qed
\end{remark}

\subsection{Lower semi-continuity and mean value properties
of extremal approximately convex functions}\label{sec:lsc}

Let $K\subset \R^n$ be a compact convex set and let $V$ be the set of
extreme points of $K$.  Let $\phi\cn V \to R$ be a function.  Then a
function $h\cn K \to \R^n$ has {\bi extreme values equal to\/} $\phi$
iff $h\big|_V=\phi$. (The terminology is a variant on that used in
partial differential equations where the boundary values of a function
are often prescribed.) Likewise if $f,g\cn K\to \R$ are two functions
then $f$ and $g$ {\bi have the same extreme values\/} iff they agree
on $V$.  If $\phi\cn V\to \R$, let ${\mathcal B}(K,\phi)$ be the
set of bounded approximately convex functions $h$ so that $h\big|_V\le
\phi$ on $V$.  Then the {\bi extremal approximately convex function\/}
with  extreme values equal to $\phi$ is
\begin{equation}\label{E-K-def}
E_{K,\phi}(x)=\sup_{h\in {\mathcal B}(K,\phi)}h(x).
\end{equation}
This is the pointwise largest approximately convex function with
extreme values $\le \phi$ on $V$. While in general we may have
$E_{K,\phi}(v)< \phi(v)$ for some $v\in V$, we will show that if
$\phi$ is uniformly continuous on $V$ (which will always be the case
if $V$ is finite) then $E_{K,\phi}\big|_V=\phi$ and that $E_{K,\phi}$
is lower semi-continuous on $K$.

Let $K\subset \R^n$ be a compact set with extreme points $V$.  Then
for any function $h\cn K\to \R$ which is bounded above define $Sh$ by
$$
Sh(x)=\begin{cases} h(x),& x\in V;\\
{\displaystyle	\inf\Bigl\{\dfrac{h(y)+h(z)}{2}+1}\ :\ 
	\dfrac{y+z}{2}=x\Bigr\},& x\in K\setminus
	V.
	\end{cases}
$$
This operator is closely related to approximately convex functions as
\begin{equation}\label{f-le-Sf}
f\le Sf\quad \iff \quad \text{$f$ is approximately convex on $K$.}
\end{equation}
Despite being nonlinear $S$ is somewhat like a mean value operator.  We
make this more precise by proving a maximum principle for the
equation $Sf=f$.

\mnote{\ralph I have tried to get all the inequalities right here.
Please check to see that I have done it correctly.}
\begin{thm}\label{S-mean}
Let $K\subset \R^n$ be a compact convex set with extreme points $V$.
Let $f,F\cn K\to \R$ be bounded functions so that  $Sf\le f$ and $F$
is approximately convex (that is $SF \ge F$).  Let
\begin{equation}\label{L-def}
L(x)=\min\{f(x),\liminf_{y\to x}f(y)\}
\end{equation}
be the lower semi-continuous envelope of $f$.  Then
\begin{equation}\label{f-sup}
\sup_{x\in K}(F(x)-f(x))=\sup_{v\in V}(F(v)-f(v))
\end{equation}
and
\begin{equation}\label{L-sup}
\sup_{x\in K}(F(x)-L(x))=\sup_{v\in V}(F(v)-L(v)).
\end{equation}
\end{thm}

\begin{proof}
We will prove~(\ref{L-sup}), the proof of~(\ref{f-sup}) being similar
(and a little easier).  The inequality $Sf\le f$ implies that if
$x\notin V$ then
\begin{equation}\label{f-inf}
f(x)\ge \inf\Bigl\{\frac{f(y)+f(z)}{2}+1\,:\, \frac{y+z}{2}=x\Bigl\}.
\end{equation}
As $f$ and $F$ are bounded we may assume (after possibly adding
positive constants to $f$ and $F$) that $0 \le f\le F\le M$ for some
positive constant $M$.  This implies $0\le L\le F$.  Set
$\omega(x):=F(x)-L(x)$ and $\delta:=\sup_{x\in K}\omega(x)$.  Then we
wish to show $\sup_{v\in V}\omega(v)=\delta$.  If $\delta=0$ then
$L\equiv F$ and there is nothing to prove.  So assume $\delta>0$.
Choose a positive integer~$N$ so that $N>M$.  Let $0<\e<1$ and choose
$w_0$ to be a point so that $\omega(w_0)>(1-\e2^{-N})\delta$. Suppose
that $w_0\notin V$ for sufficiently small $\e>0$ (otherwise the
desired conclusion follows as $\e\to 0$). From the definition of $L$
there is a sequence $\la x_k\ra_{k=1}^\infty$ such that $x_k\to w_0$
and $f(x_k)\to L(w_0)$. By equation~(\ref{f-inf}) there are sequences
$\la y_k\ra_{k=1}^\infty$ and $\la z_k\ra_{k=1}^\infty$ such that
$x_k=(y_k+z_k)/2$ and a real number $C\ge 0$ such that
\begin{equation}\label{f-lim}
f(x_k)-\Bigl(\frac{f(y_k)+f(z_k)}{2}+1\Bigr) \to C\ge 0
\end{equation}
By passing to a subsequence we may assume that $y_k\to y$, $z_k\to z$,
$f(y_k)\to A$, and $f(z_k)\to B$ for some $y,z\in K$ and $A,B\in \R$.
Clearly $w_0=(y+z)/2$ and (using the definition of $L$) $L(y)\le A$,
$L(z)\le B$. Then~(\ref{f-lim}) yields
\begin{equation}\label{L(w)}
L(w_0)=\frac{A+B}{2}+1+C\ge \frac{L(y)+L(z)}{2}+1
\end{equation}
and so
\begin{equation}\label{L(w)-lower}
F(w_0)=L(w_0)+\omega(w_0)\ge \frac{L(y)+L(z)}{2}+1+\omega(w_0).
\end{equation}
But since $F$ is approximately convex 
\begin{equation}\label{L(w)-upper}
F(w_0)\le \frac{F(y)+F(z)}{2}+1
=\frac{L(y)+L(z)}{2}+1+\frac{\omega(y)+\omega(z)}{2}.
\end{equation}
Combining~(\ref{L(w)-upper}) and~(\ref{L(w)-lower}) yields
\begin{equation}\label{omega-bd}
\frac{\omega(y)+\omega(z)}{2}\ge \omega(w_0).
\end{equation}
Since $\delta=\sup_{x\in K}\omega(x)$ and $\omega(w_0)\ge
(1-\e2^{-N})\delta$, (\ref{omega-bd}) implies
$$
\min\{\omega(y),\omega(z)\}\ge 2\omega(w_0)-\delta
\ge  (2(1-\e2^{-N})-1)\delta=(1-\e2^{-(N-1)})\delta.
$$ 
From~(\ref{L(w)}) have $\min\{L(y),L(z)\}\le
L(w_0)-1$.  Without loss of generality we may assume that $L(y)\le
L(w_0)-1$.  Let $w_1=y$.  Then $L(w_1)\le L(w_0)-1$ and
$\omega(w_1)\ge(1-\e2^{-(N-1)})$.

If $w_1\notin V$ then we can repeat this argument (with $N$ replaced
by $N-1$) and get a $w_2\in K$ with $L(w_2)\le L(w_1)-1$ and
$\omega(w_1)\ge (1-\e2^{N-2})\delta$.  We continue in this manner to
get a finite sequence $w_0, w_1\cd w_m$ with $m<N$ so that for $1\le
k\le m-1$ we have $L(w_k)\le L(w_{k-1})-1$, $\omega(w_k)\ge
(1-\e2^{N-k})\delta$, and $w_k\notin V$. (Note this can not
continue for $k\ge N$ as that would imply $L(w_N)\le L(w_0)-N\le
M-N<0$ contradicting $L\ge 0$.  Thus $w_m\in V$ for some $m<N$.) At
the last step $w_m\in V$ and $\omega(w_m)\ge (1-\e 2^{-(N-m)})\delta$.
Therefore $\sup_{v\in V}\omega(v)\ge (1-\e 2^{-(N-m)})\delta$.  Letting
$\e\searrow 0$ yields $\sup_{v\in V}\omega(v)\ge \delta$.  But
$\sup_{v\in V}\omega(v)\le \delta$ is clear.  Thus $\sup_{v\in
V}\omega(v)= \delta$ as required.
\end{proof}

\begin{prop}\label{S-basics}
Let $K\subset \R^n$ be a compact convex set with extreme points $V$
and let $h\cn K \to \R^n$ be a bounded approximately convex function
on $K$.  Then $h(x)\le Sh(x)$, the functions $h$ and $Sh$ have the
same extreme values, $Sh$ is approximately convex, and if $h$ is lower
semi-continuous as a function on $K$ at points of $V$ then the same is
true of $Sh$.
\end{prop}

\begin{proof}
If $x=(y+z)/2$ then as $h$ is approximately convex $h(x)\le
(h(y)+h(z))/2+1$ and taking the infimum yields $h(x)\le Sh(x)$. That
$h$ and $Sh$ have the same extreme values is clear.  Using the
definition of $Sh$ and the inequality $h\le Sh$ we have
$$
Sh\Bigl(\frac{y+z}{2}\Bigr)\le 
	\frac{h(y)+h(z)}{2}+1\le \frac{Sh(y)+Sh(z)}{2}+1,
$$
which shows $Sh$ is approximately convex.  Finally if $h$ is lower
semi-continuous at points of $V$ then for $x\in V$ we have
$\liminf_{y\to x}Sh(y)\ge \liminf_{y\to x}h(y)\ge h(x)=Sh(x)$. This
shows $Sh$ is lower semi-continuous at $x$ and completes the proof.
\end{proof}

We now characterize the extremal functions $E_{K,\phi}$ as the unique
bounded  solutions to the equation $Sf=f$ with extreme values $\phi$.

\begin{thm}\label{Sf=E}
Let $K$ be a convex set with extreme points $V$ and $f\cn K\to \R$ a
bounded function so that $Sf=f$.  Let $\phi:=f\big|_V$ be the extreme
values of $f$ and let $E_{K,\phi}$ be the extremal approximately convex
function with extreme values $\phi$.  Then $f=E_{K,\phi}$.
\end{thm}

\begin{proof}
The equality $Sf=f$ implies $f$ is approximately
convex~(cf.~(\ref{f-le-Sf})).  Then the extremal property of
$E_{K,\phi}$ implies $f\le E_{K,\phi}$.  Let $F=E_{K,\phi}$ in
Theorem~\ref{S-mean} and using that $f$ and $E_{K,\phi}$ agree on $V$
we can use equation~(\ref{f-sup}) to conclude $f=E_{K,\phi}$.
\end{proof}

The following is an elementary variant on Corollary~17.2.1
in~\cite{Rockafellar:convex}.  We include a short proof for
completeness.

\begin{prop}\label{convex-bdry}
Assume $K\subset \R^n$ is a compact convex set and $V$ the set of
extreme points of $K$. Let $\phi\cn V\to \R$ be uniformly continuous.
Then there exists a lower semi-continuous convex function $h\cn K\to
\R$ so that $h\big|_V=\phi$.  Moreover we can choose $h$ so that
$\inf_{x\in K}h(x)=\inf_{v\in V}\phi(v)$ and $\sup_{x\in
K}h(x)=\sup_{v\in V}\phi(v)$.
\end{prop}

\begin{proof}
Let $\ol{V}$ be the closure of $V$.  As $\phi\cn V\to \R$ is uniformly
continuous it has a unique continuous extension $\ol{\phi}\cn
\ol{V}\to \R$.  Let Let $G_{\ol{\phi}}:=\{(x,\ol{\phi}(x)): x\in
\ol{V}\}\subset K\times \R$ be the graph of $\ol{\phi}$.  As the set
$\ol{V}$ is a compact and $\ol{\phi}$ is continuous the set
$G_{\ol{\phi}}$ is also compact.  Therefore the convex hull
$\conv(G_{\ol{\phi}})$ is compact. Let $A:=\inf_{v\in
V}\phi(v)=\min_{x\in \ol{V}}\ol{\phi}(x)$ and $B:=\sup_{v\in
V}\phi(v)=\max_{x\in \ol{V}}\ol{\phi}(x)$.  Then
$\conv(G_{\ol{\phi}})\subseteq K\times [A,B]$.  Moreover, as $K$ is
the convex hull of its set of extreme points $V$,  if $x\in K$
then there is $y\in [A,B]$ so that $(x,y)\in \conv(G_\phi)$.  Define
$h$ by
$$
h(x):=\min\{y: (x,y)\in \conv(G_{\ol{\phi}})\}.
$$
It is clear from this definition that $h$ is convex and has the same
supremum and infimum as $\phi$.  We now show that $h$ is lower
semi-continuous.  Let $a\in K$ and let $A:=\liminf_{x\to
a}f(x)$. Choose a sequence $\la x_\ell\ra_{\ell=1}^\infty$ so that
$x_\ell\to a$ and $h(x_\ell)\to A$.  Then as $\conv(G_{\ol{\phi}})$ is
compact (and thus closed) the limit $\lim_{\ell\to
\infty}(x_\ell,h(x_\ell))=(a,A)\in \conv(G_{\ol{\phi}})$.  The definition of
$h$ then implies $h(a)\le A=\liminf_{x\to a}h(x)$.  Thus $h$ is lower
semi-continuous at $a$ for every $a\in A$.

Finally let $v\in V$.  Then as $(v,h(v))\in \conv(G_{\ol{\phi}})$
there exists $(\alpha_0\cd \alpha_{n+1})\in \Delta_{n+1}$ and $v_0\cd
v_{n+1}\in \ol{V}$ so that
$(v,h(v))=\sum_{k=0}^n\alpha_k(v_k,\ol{\phi}(v_k))$.  But $v$ is an
extreme point of $K$, which implies that $v_k=v$ for all $k$ and
therefore $h(v)=\ol{\phi}(v)=\phi(v)$.
\end{proof}

\begin{thm}\label{extreme-lsc}
Let $K\subseteq \R^n$ be a compact convex set with extreme points $V$.
Assume that $\phi\cn V\to \R$ is uniformly continuous.  Then the
extremal approximately convex function $E_{K,\phi}$ satisfies
$E_{K,\phi}\big|_V=\phi$ and is lower semi-continuous on $K$.
\end{thm}

\begin{proof}
By Proposition~\ref{convex-bdry} there exists a lower semi-continuous
convex function $h\cn K\to \R$ with extreme values $\phi$.  As $h$ is
convex it is {\it a fortiori\/} approximately convex. $h$ 
approximately convex (so that $h\le E_{K,\phi}$) we have for $v\in V$
that $\phi(v)=h(v) \le E_{K,\phi}(v)\le \phi(v)$, and so $E_{K,\phi}$
has $\phi$ as extreme values.  As $h\le E_{K,\phi}$ and $h$ is lower
semi-continuous, the function $E_{K,\phi}$ will be lower
semi-continuous at all points $x$ where $E_{K,\phi}(x)=h(x)$.  In
particular, $E_{K,\phi}$ will be lower semi-continuous at all points of
$V$.  Finally as $SE_{K,\phi}\ge E_{K,\phi}$ (cf.~\ref{S-basics}) the
extremal property of $E_{K,\phi}$ implies $SE_{K,\phi}=E_{K,\phi}$.
Now in Theorem~\ref{S-mean} let $f=F=E_{K,\phi}$ and let $L$ be the
lower semi-continuous envelope of $f=E_{K,\phi}$ as given
by~(\ref{L-def}).  Then as $E_{K,\phi}$ is lower semi-continuous at
points of $V$ we have that $E_{K,\phi}(v)=L(v)$ for all $v\in V$.
Therefore~(\ref{L-sup}) implies that $E_{K,\phi}=L$ on $K$, so that
$E_{K,\phi}$ is lower semi-continuous as claimed.
\end{proof}

\begin{remark}\label{E=S-Limit}
Let $K\subset \R^n$ be a convex set with extreme points $V$.  Let
$h\cn K\to \R$ be a bounded approximately convex function and let
$\phi\cn V\to \R$ be the extreme values of $h$, that is
$\phi:=h\big|_V$.  Then there is a bounded function $f\cn K\to \R$
such that $f\big|_V=\phi$ for which the inequality $Sf\le f$ holds
pointwise on $K$. (Such a function exists as is seen by letting
$f=E_{K,\phi}$.  On the simplex $\Delta_n$ with $\phi=0$ the function
$f(x)=k$ for $x$ in the interior of a $k$-dimensional face is
an example of such a function.)  Then define two sequences $\la
h_k\ra_{k=0}^\infty$ and $\la f_k\ra_{k=0}^\infty$ of functions on $K$
by
$$
h_0=h,\  h_{k+1}=Sh_k,\quad f_0=f,\  f_{k+1}=Sf_k.
$$
Then it can be shown that $f_{k+1}\le f_k$, $h_{k+1}\ge h_k$, and that
each $h_k$ is approximately convex.  (The statements about $h_k$
follow from Proposition~\ref{S-basics}.)  Also all the $h_k$'s and
$f_k$'s have $\phi$ as extreme values.  Therefore both sequences have
pointwise limits $h_\infty=\lim_{k\to \infty}h_k$ and
$f_\infty=\lim_{k\to\infty}f_k$.  These both have $\phi$ as extreme
values, $Sh_\infty=h_\infty$, and $Sf_\infty=f_\infty$.  Therefore by
Theorem~\ref{Sf=E} we have $h_\infty=f_\infty=E_{K,\phi}$.  This gives
a method for finding $E_{K,\phi}$ as the limit of two more or less
constructively defined sequences.  Also note that for each $k$ we have
the inequalities
$$
h_k\le E_{K,\phi}\le f_k.
$$
Thus we have explicit upper and lower bounds for $E_{K,\phi}$.\qed
\end{remark}

\subsection{The extremal approximately sub-affine function $H(x)$}
\label{sec:H}

A function $f\cn [0,1]\to \R$ is {\bi approximately sub-affine\/} iff
\begin{equation}\label{sub-lin}
f\Bigl(\frac{x+y}{2}\Bigr)\le \frac{f(x)+f(y)}{2}+\frac{x+y}{2}.
\end{equation}
As in example~\ref{sub-convex} below approximately sub-affine
functions can be used to construct approximately convex functions on a
simplex.  As a first step in explicitly describing the extremal
approximately convex function on a simplex we describe the extremal
approximately convex function on the unit interval.

\begin{figure}[hb]
\centering
\mbox{\psfig{file=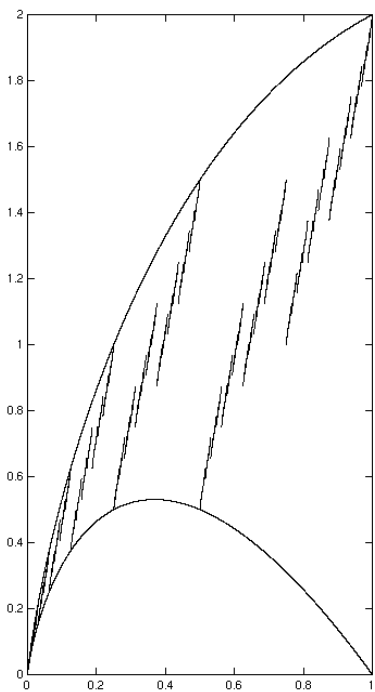,width=2in}}
\caption[]{Graphs of $y=H(x)$, $y=x\log_2(x)$, and
$y=2x+\log_2(x)$} for $0\le x\le 1$.
\label{H-graph}
\end{figure}

Let $\Nat=\{0,1,2,\dots\}$ be the natural numbers and let $\dy$ be
the {\bi dyadic rational numbers\/} in $[0,1]$.  That is
$$
\dy :=\Bigr\{\frac{m}{2^n}: m, n\in \Nat \text{ and } 0\le m \le 2^n\Bigl\}.
$$
(These  play a considerable r\^ole in what follows.)  The numbers in
$[0,1]\setminus \dy$ will be called the {\bi dyadic irrationals\/}.
Every dyadic irrational $x$ has a unique binary expansion
$x=\sum_{i=0}^\infty {x_i}/{2^i}$ with $x_i\in\{0,1\}$.  If $x\in
\dy$ then there are two binary expansions: the finite expansion 
$ x=\sum_{i=0}^N {x_i}/{2^i} $
and, if $x_N=1$, there is also the infinite expansion 
$x=\sum_{i=0}^{N-1} {x_i}/{2^i}+\sum_{i=N+1}^\infty
{1}/{2^i}$.  Unless stated otherwise we will always use the finite
expansion for an element of $\dy$, even when we write
$x=\sum_{i=0}^\infty x_i/2^i$ for notational uniformity.  With this
understood, define $H\cn [0,1]\to \R$ by
\begin{equation}\label{H-def}
H(x):=\sum_{i=0}^\infty i\frac{x_i}{2^i} \quad \text{where}\quad	
	 x=\sum_{i=0}^\infty\frac{x_i}{2^i}.
\end{equation}
For motivation see Remark~\ref{dy-H-motivation}.  A graph of $H$ is
shown in Figure~\ref{H-graph}.

We now derive another representation of $H$. 
Let $r\cn \R\to \R$ be defined by
$$
r(x):=\begin{cases} 0,& 0\le x< 1; \\ 1, & 1\le x<2,\end{cases}
$$
and  extend to $\R$ by periodicity: $r(x+2)=r(x)$.
If $0\le x<1$ and $x$ has binary expansion $x=\sum_{i=1}^\infty
x_i/2^i$, where $x_i\in \{0,1\}$, then it is not hard to see that
$x_i=r(2^{i}x)$ (if $x$ is a dyadic rational we check to see this does
give the finite expansion).  It follows  for $0\le x <1$ that
$x=\sum_{i=1}^\infty r(2^ix)/2^i$.  More generally if we let
$\fp{x}=x-[x]$ be the fractional part of $x$ then as both $\fp{x}$
and $\sum_{i=1}^\infty r(2^ix)/2^i$ are periodic with period $1$ and
$\fp{x}=x$ for $0\le x<1$ we have
\begin{equation}\label{frac-series}
\fp x=\sum_{i=1}^\infty \frac{r(2^ix)}{2^i}.
\end{equation}
If $H$ is extended to $\R$ to be periodic, $H(x+1)=H(x)$, (this is
possible as $H(0)=H(1)=0$) then the definition of $H$ becomes
\begin{equation}\label{H-def-rad}
H(x)=\sum_{i=1}^\infty  i\frac{r(2^ix)}{2^i}.
\end{equation}

\begin{prop}\label{H-series}
Let the function $H$ be extended from $[0,1)$ to $\R$ so that $H$ is
periodic: $H(x+1)=H(x)$.  Then $H$ satisfies the functional equation
\begin{equation}\label{H-fcn-eqn}
H(x)=\fp x +\frac{1}{2}H(2x)
\end{equation}
and thus $H$ has the series representation
\begin{equation}\label{H-frac-series}
H(x)=\sum_{k=0}^\infty \frac{\fp{2^kx}}{2^k}.  
\end{equation}
This implies  $H$ is lower semi-continuous,  continuous at all
points of $[0,1]\setminus \dy$ and right continuous at all points. 
\end{prop}

\begin{proof} This is a calculation based on the two
series~(\ref{H-def-rad}) and~(\ref{frac-series}).
\begin{align*}
H(x)&=\sum_{i=1}^\infty i\frac{r(2^ix)}{2^i}
=\fp x + \sum_{i=1}^\infty (i-1)\frac{r(2^ix)}{2^i}\\
&=\fp x + \frac12 \sum_{i=2}^\infty (i-1)\frac{r(2^{i-1}2x)}{2^{i-1}}\\
&=\fp x + \frac12 \sum_{j=1}^\infty j\frac{r(2^{j}2x)}{2^{j}}
=\fp x + \frac12 H(2x). 
\end{align*}
To prove the series representation~(\ref{H-frac-series}) for $H(x)$
observe that an induction using the functional equation~(\ref{H-fcn-eqn})
yields
$$
H(x)=\sum_{k=0}^m\frac{\fp{2^kx}}{2^k} +\frac{1}{2^{m+1}}H(2^{m+1}x)
$$
and as $0\le H(x)\le \sum_{i=1}^\infty i/2^i=2$ the series converges
uniformly to $H(x)$.  The functions $x \mapsto {\fp{2^kx}}/{2^k}$ are
lower semi-continuous and right continuous and hence so are the
partial sums $S_n(x)=\sum_{k=0}^n{\fp{2^kx}}/{2^k}$.  Thus $H$ is the
uniform limit of lower semi-continuous and right continuous functions
and therefore is lower semi-continuous and right continuous.  Finally
the functions ${\fp{2^kx}}/{2^k}$ are continuous at all points of
$[0,1]\setminus \dy$.  As the series converges uniformly this implies
that the sum $H$ is also continuous at these points.
\end{proof}

\begin{remark}\label{self-cong}
The graph of $H(x)$ has an interesting ``self-congruence'' property.
The series~(\ref{H-frac-series}) for $H(x)$ implies  for $m$ a
positive integer that
$$
H\Bigl(x+\frac1{2^m}\Bigr)
	 =\sum_{k=0}^{m-1}\frac1{2^k}(\fp{2^kx+2^{k-m}}-\fp{2^kx})+H(x)
	 =P_m(x)+H(x)
$$ 
where this defines $P_m(x)$.  It is not hard to check that the
functions $(\fp{2^kx+2^{k-m}}-\fp{2^kx})/2^k$ are all constant on
intervals $[i/2^m,(i+1)/2^m)$ and so the same will be true for
$P_m(x)$.  This implies for any $i$ and $j$ that the graph of the
restriction $H\big|_{[i/2^m,(i+1)/2^m)}$ is a translation of the graph
of $H\big|_{[j/2^m,(j+1)/2^m)}$.  So informally and somewhat
imprecisely ``the graph of $H$ is locally self congruent at all the
scales $1/2^m$''. If $F$ is the closure of the graph of
$H\big|_{[0,1)}$ then this, and some calculation, can be used to show
$F$ can be covered by $2^m$ closed sets of diameter $\le 4
m2^{-m}$. Thus for any $\delta>0$ the Hausdorff $\delta$-dimensional
measure of $F$ is $\le 2^m(4m2^{-m})^\delta$ and when $\delta >1$ we
have $2^m(4m2^{-m})^\delta\to 0$ as $m\to \infty$.  Therefore the
Hausdorff dimension of $F$ is $\le1$.  But as $F$ projects onto the
interval $[0,1]$ its Hausdorff dimension is $\ge 1$.  Thus $F$ has
Hausdorff dimension one.  (With a little more work it can be shown the
one dimensional Hausdorff measure of $F$ is infinite.)  However $F$ is
compact, separable, totally disconnected and has no isolated points.
Thus $F$ is homeomorphic to the Cantor set and therefore of
topological dimension zero.  Whence the closure of the graph of $F$ is
a ``fractal'' in the sense that its geometric dimension is greater
than its topological dimension.\qed
\end{remark}

\begin{prop}\label{H-sub-affine}
The function $H$ is approximately sub-affine:
$$
H\Bigl(\frac{x+y}2\Bigr)\le \frac{H(x)+H(y)}{2}+\frac{x+y}{2}\quad \text{for}\quad
x,y\in [0,1].
$$
\end{prop}

\begin{lemma} \label{II1}
If $x=\sum_{i=0}^N{l_i}/{2^i}\in \dy$ with each $l_i$ a
nonnegative integer, then
$$
H(x)\le \sum_{i=0}^Ni\frac{l_i}{2^i}
$$
with equality if and only if each $l_i\in \{0,1\}$.
\end{lemma}

\begin{proof}
If $(l_0,l_1\cd l_N)$ is a finite sequence with
$\sum_{i=0}^N{l_i}/{2^i}\le 1$ we let $\lambda(l_0\cd
l_N):=\sum_{i=0}^Nl_i$.  The proof is by induction on
$m=\lambda(l_0\cd l_N)$.  If $m=0$ then each $l_i=0$ and $x=H(x)=0$
and the result is trivial.  Now assume the inequality holds for all
$(l_0,l_1\cd l_N)$ with $\lambda(l_0\cd l_N)<m$.  Let $k$ be the least
integer such that $l_k\ge2$ (if all $l_k\in\{0,1\}$ there is nothing to
prove).  Note $k\ne0$ as $2/2^0=2$. Then
\begin{align*}
x&=\sum_{i=0}^{k-2}\frac{l_i}{2^i} +\frac{l_{k-1}}{2^{k-1}} 
	 +\frac{l_{k}}{2^{k}}+ \sum_{i=k+1}^{N}\frac{l_i}{2^i}\\
&=\sum_{i=0}^{k-2}\frac{l_i}{2^i} +\frac{l_{k-1}+1}{2^{k-1}} 
	 +\frac{l_{k}-2}{2^{k}}+ \sum_{i=k+1}^{N}\frac{l_i}{2^i}
=\sum_{i=0}^N\frac{r_i}{2^i}
\end{align*}
where the last line defines the $r_i$ implicitly.  Then
$$
\lambda(r_0\cd r_N)=\lambda(l_0\cd l_{k-1}+1,l_k-2\cd l_N)=
\lambda(l_0\cd l_N)-1=m-1.
$$
Thus the induction hypothesis gives
\begin{align*}
\sum_{i=0}^N i\frac{l_i}{2^i}
&=\sum_{i=0}^N i \frac{r_i}{2^i}+k\frac{2}{2^k}-(k-1)\frac{1}{2^{k-1}}\\
&=\sum_{i=0}^N i \frac{r_i}{2^i} +\frac{k-(k-1)}{2^{k-1}}
> \sum_{i=0}^N i \frac{r_i}{2^i} \ge H(x).
\end{align*}
This gives $H(x)<\sum_{i=0}^N {l_i}/{2^i}$ unless 
$l_i\in \{0,1\}$ for all $i$.  This completes the proof.
\end{proof}

\begin{proof}[Proof of Proposition~\ref{H-sub-affine}]
First consider the case $x,y\in \dy$
so that $x=\sum_{i=0}^Nx_i/2^i$, $y=\sum_{i=0}^Ny_i/2^i$.
Then by Lemma~\ref{II1}
\begin{align*}
H\Bigl(\frac{x+y}{2}\Bigr)&=H\bigg(\sum_{i=0}^N\frac{x_i+y_i}{2^{i+1}}\biggr)
\le \sum_{i=0}^N(i+1)\frac{x_i+y_i}{2^{i+1}}\\
&=\frac12\biggl(\sum_{i=0}^Ni\frac{x_i}{2^i}+\sum_{i=1}^Ni\frac{y_i}{2^i}\biggr)
+\frac12\biggl(\sum_{i=0}^N\frac{x_i}{2^i}+\sum_{i=1}^N\frac{y_i}{2^i}\biggr)\\
&=\frac{H(x)+H(y)}{2}+\frac{x+y}2
\end{align*}
If $x$ is a dyadic irrational and $y\in\dy$ then we use that by
Proposition~\ref{H-series} the function $H$ is lower semi-continuous
on $\R$ and continuous at $x$.  Let $x(r)\in\dy$ so that
$\lim_{r\to\infty}{x(r)}=x$ and so by continuity
$\lim_{r\to\infty}H(x(r))=H(x)$.  Thus
\begin{align*}
H\Bigl(\frac{x+y}{2}\Bigr)&
	 \le \liminf_{r\to \infty}H\Bigl(\frac{x(r)+y}{2}\Bigr)
\le \lim_{r\to\infty}\(\frac{H(x(r))+H(y)}{2}+\frac{x(r)+y}{2}\)\\
&=\frac{H(x)+H(y)}{2}+\frac{x+y}{2}
\end{align*}
The case where both $x$ and $y$ are dyadic irrationals is handled
similarly.
\end{proof}

\begin{prop}\label{H-biggest}
Suppose $f$ is a lower semi-continuous approximately sub-affine
function defined on $[0,1]$ such that $f(0)=0$.  Then $f(x)\le H(x)+f(1)x$
for all $x\in [0,1]$.
\end{prop}

First some preliminaries.  If $x=\sum_{j=1}^Nx_j/2^j\in \dy$ 
define the {\bi dyadic support\/} of $x$ to be $\{j\in \Nat: x_j=1\}$
and denote it by $\supp x$.

\begin{lemma}\label{H-exact}
If $x,y\in \dy$ and $(\supp x)\cap (\supp y)=\nothing$ then 
$$
H\Big(\frac{x+y}{2}\Bigr)=\frac{H(x)+H(y)}{2}+\frac{x+y}{2}.
$$
\end{lemma}

\begin{remark}\label{dy-H-motivation}
This lemma motivated the definition of $H$.  As the proof of
Proposition~\ref{H-biggest} makes clear this is the property which
implies $H$ is the largest lower semi-continuous approximately sub-affine
function on $[0,1]$.   It also allows one to compute the values of
$H$ on $\dy$ leading to the formula~(\ref{H-def}).\qed
\end{remark}

\begin{proof} The condition on the dyadic supports implies
that the binary expansion of $x+y$ can be computed by just adding the
digits without ``carrying''.  Thus for sufficiently large $N$
\begin{align*}
H\Bigl(\frac{x+y}{2}\Bigr)&=\sum_{j=0}^N(j+1)\frac{x_j+y_y}{2^{j+1}}\\
&=\frac12\biggl(\sum_{j=0}^N(j+1)\frac{x_j}{2^{j}}
	 +\sum_{j=0}^N(j+1)\frac{y_j}{2^{j}}\biggr)+\frac{x+y}{2}\\
&=\frac{H(x)+H(y)}{2}+\frac{x+y}{2}.
\end{align*}
\vskip-15pt
\end{proof}

\begin{proof}[Proof of Proposition~\ref{H-biggest}] 
If $f(x)$ is replaced by $\phi(x):=f(x)-f(1)x$ then $\phi$ will also
be approximately sub-affine and $\phi(0)=\phi(1)=0=H(0)=H(1)$.  We now
show by induction on $k$ that if $x=m/2^k\in \dy$ then $\phi(x)\le
H(x)$.  The base case of $k=0$ holds.  Now assume that $x=m/2^k$ and
that the result is true when the denominator of the fraction is a
smaller power of~$2$.  We may assume that $m$ is odd.  If $x\le 1/2$
let $y=2x=m/2^{k-1}$.  Then $x=(0+y)/2$, $\phi(y)\le H(y)$ and
$\supp(0)\cap \supp(y)=\nothing$.  Therefore
\begin{align*}
\phi(x)&=\phi\Bigr(\frac{0+y}{2}\Bigr)\le
\frac{\phi(0)+\phi(y)}{2}+\frac{0+y}{2}\\
&\le \frac{H(0)+H(y)}{2}+\frac{0+y}{2}=H\Bigl(\frac{0+y}{2}\Bigr)=H(x).
\end{align*}
If $1/2<x<1$ then let $y=2x-1$ so that $x=(y+1)/2$.  Then as the
dyadic supports of $y$ and $1$ are disjoint, a calculation like the one
just done shows $\phi(x)\le H(x)$.  Thus $\phi(x)\le H(x)$ for
all $x\in \dy$.  For any other 
$x\in [0,1]\setminus \dy$ choose $x_k\in \dy$ with
$x_k\to x$.  By Proposition~\ref{H-series} $H$  is continuous at
$x$.  Therefore the lower semi-continuity of $\phi$ implies
$$
\phi(x)\le \liminf_{k\to \infty}\phi(x_k)\le \lim_{k\to\infty}H(x_k)=H(x).
$$
Finally $\phi(x)\le H(x)$ is equivalent to the required inequality for
$f$.
\end{proof}

\begin{prop}\label{H-bounds}
The inequalities
\begin{equation}\label{<H<}
x\log_2(1/x)\le H(x) \le 2x+\log_2(1/x)
\end{equation}
hold for $0\le x\le 1$ (cf.~Figure~\ref{H-graph}).
\end{prop}

\begin{lemma}\label{x*log(x)-ineq}
Let $\phi(x):=x\log_2(1/x)=-x\ln(x)/\ln(2)$.  Then for $0\le t \le 1$
and $x,y\in [0,1]$
$$
0\le \phi((1-t)x +t x) - t\phi(x) -(1-t)\phi(x) \le \phi(t)x+\phi(1-t)y.
$$
As $\phi(1/2)=1/2$ this implies $\phi$ is approximately sub-affine on
$[0,1]$.
\end{lemma}

\begin{proof}  The left hand inequality follows from the concavity of
$\phi$. To prove the right hand inequality we first assume $0<x\le
y\le 1$.  For fixed $t$ and $y$ let
$$
F(x):=\phi((1-t)x+ty)-(1-t)\phi(x)-t\phi(y).
$$
Then
$$
F'(x)=\frac{(1-t)}{\ln(2)}(\ln(x)-\ln((1-t)x+ty))\le 0.
$$
Therefore $F$ is monotone decreasing and so the maximum of $F(x)$ on
$[0,y]$ occurs when $x=0$.  But
$$
F(0)=\phi(ty)-t\phi(y)=\frac{-(ty\ln(ty)-ty\ln(y))}{\ln(2)}
	 =\frac{-t\ln(t)}{\ln(2)}\,y=\phi(t)y.
$$
So for all $0\le x\le y$ and $0\le t\le 1$
$$
\phi((1-t)x+ty)-t\phi(x)-(1-t)\phi(y)\le \phi(t)y\le \phi(t)y+\phi(1-t)x
$$
(for the last step note that  $\phi(1-t)x\ge 0$). A
similar argument works in the case $y\le x$ (or replace $t$ by $(1-t)$
in what has been shown).
\end{proof}

\begin{proof}[Proof of Proposition~\ref{H-bounds}]
As the function $\phi(x)=x\log_2(1/x)$ is approximately sub-affine,
vanishes at the endpoints of $[0,1]$ and is continuous the lower bound
of (\ref{<H<}) follows from Proposition~\ref{H-biggest}.  To prove the
upper bound we use the series~(\ref{H-frac-series}).  Let $0<x<1$.
There exists a unique nonnegative integer~$m$ so that $2^mx<1\le
2^{m+1}x$ (i.e. $1/2^{m+1}\le x< 1/2^m$).  Then for $0\le k\le m$ we
have $\fp{2^kx}=2^kx$, and thus
$$
H(x)=\sum_{k=0}^\infty\frac{1}{2^k}\fp{2^kx}\le
(m+1)x+\sum_{k=m+1}^\infty \frac{1}{2^k}=(m+1)x+\frac{1}{2^m}.
$$
So to complete the proof it is enough to show 
$$
\psi(x):=2x+x\log_2(1/x)-\Bigl((m+1)x+\frac{1}{2^m}\Bigr)
=\frac{-x\ln(x)}{\ln(2)}-(m-1)x+\frac{1}{2^m}
$$
satisfies $\psi(x)\ge 0$ for $x\in [1/2^{m+1},1/2^m]$.  But
$\psi(1/2^{m+1})=\psi(1/2^m)=0$ and $\psi''(x)=-1/(x\ln(2))<0$.  So
$\psi$ is concave on $[1/2^{m+1},1/2^m]$ and vanishes at the endpoints
which implies $\psi\ge0$ on the interval.
\end{proof}

\subsection{The extremal approximately convex function $E(x)$ on a simplex}
\label{sec:E}

Let $e_0\cd e_n$ be the standard basis of $\R^{n+1}$.   Then the standard
simplex is, as usual, $\Delta_n=\conv\{e_0\cd e_n\}$.  We will often
write points of $\Delta_n$ in terms of their affine coordinates
$(x_0\cd x_n)$ where $x_k\ge 0$ and $\sum_{k=0}^nx_k=1$.  This
corresponds to $\sum_{k=0}^nx_ke_k$.  
Define a
function $E$ on $\Delta_n$ as follows:
\begin{equation}\label{E-defn}
E\biggl(\sum_{k=0}^nx_ke_k\biggr)=E(x_0\cd x_n):=\sum_{k=0}^nH(x_k).
\end{equation}

\mnote{\ralph FIND GOOD REFERENCE.}
\begin{remark}\label{E-motivation}
If $\mu$ is a finite measure space and $\mathcal A$ is a finite
algebra of measurable sets with atoms $A_0, A_1\cd A_n$ the {\bi
entropy\/} of $\mathcal A$ is $ -\sum_{k=0}^n\mu(A_k)\ln\mu(A_k)$.  If
$x\in \Delta_n$ we can think of $x$ as a measure on $\{0,1,\cd
n\}$. If $\mathcal A$ is the algebra of subsets of $\{0,1,\cd n\}$
then its entropy with respect to the measure determined by $x$ is
$-\sum_{x_k}^nx_k\ln x_k$.  By Lemma~\ref{x*log(x)-ineq} the function
$x\log_2(1/x)$ is approximately sub-affine and so $H$ can be viewed as
an extremal version of $x\log_2(1/x)$.  To the extent that $H(x)$ and
$-x\ln(x)$ can be thought of as analogous functions,
$E(x)=\sum_{k=0}^nH(x_k)$ can be viewed as a ``poor man's'' version of
the entropy.  The inequalities~\ref{<H<} make this analogy somewhat
precise.  \qed
\end{remark}

The {\bi standard dyadic simplex\/} is
$$
\dy_n:=\biggl\{\sum_{k=0}^n x_ke_k: x_k\in \dy,\  \sum_{k=0}^nx_k=1\biggr\}.
$$
Like $\dy\subset [0,1]$ the set $\dy_n$ will play a large r\^ole.

\begin{prop}\label{E-basics}
The function $E$ is approximately convex and lower semi-continuous on
$\Delta_n$ with $E(e_k)=0$ for $0\le k\le n$.  The points of
continuity of $E$ are the points $x=(x_0\cd x_n)$ such that all the
coordinates $x_k$ are dyadic irrationals.  Moreover $E$ satisfies the
inequalities
$$
\sum_{k=0}^nx_k\log_2(1/x_k)\le E(x_0e_0+\cdots+x_ne_n)
	 \le 2+\sum_{k=0}^nx_k\log_2(1/x_k).
$$
\end{prop}

\begin{proof}
For $x\in \Delta_n$ the functions $x\mapsto H(x_k)$ are lower
semi-continuous by Proposition~\ref{H-series}.  Thus $E$ will also be
lower semi-continuous.  Also from Proposition the points of continuity
of $H$ are the dyadic irrationals in $[0,1]$.  This implies the
statement about the points of continuity of $E$.
As $H$ is approximately sub-affine we have 
\begin{align*}
E\Bigl(\frac{x+y}{2}\Bigr)&=\sum_{k=0}^nH\Bigl(\frac{x_k+y_k}{2}\Bigr)\\
&\le
\sum_{k=0}^n\frac{H(x_k)+H(y_k)}{2}+\sum_{k=0}^n\frac{x_k+y_k}{2}\\
&=\frac{E(x)+E(y)}{2}+1
\end{align*}
as $\sum_{k=0}^nx_k=\sum_{k=0}^ny_k=1$.  So $E$ is approximately convex as
claimed.  That $E(e_k)=0$ follows from $H(0)=H(1)=0$.  The bounds for
$E$ follow from the inequalities~(\ref{<H<}).
\end{proof}

It is possible to give an explicit formula for $E$ on the one
dimensional simplex.

\begin{prop}\label{1D-E-formula}
Let the one dimensional simplex $\Delta_1$ be identified with $[0,1]$
in the usual manner ($t$ corresponds to $(1-t)e_0+te_1$).
Then
\begin{equation}\label{easy-E}
E(t)=\begin{cases} 2,& t\notin \dy;\\ 2-\dfrac{1}{2^{l-1}},&
	 \dfrac{m}{2^l}\in \dy \text{ with $m$ odd}.\end{cases}
\end{equation}
\end{prop}

\begin{proof}  Set $\psi(t)=\fp{t}+\fp{1-t}=\fp t +\fp{-t}$.  Then
by~(\ref{H-frac-series}) 
\begin{equation}\label{E-psi}
E(t)=H(t)+H(1-t)=\sum_{k=0}^\infty \frac{\psi(2^kt)}{2^k}.
\end{equation}
But then $\psi(t)=0$ for $t\in
\mathbf{Z}$ and $\psi(t)=1$ for $t\notin \mathbf{Z}$.
So if $t\notin \dy$ we have $\psi(2^kt)=1$ for all $k$.  If
$t=m/2^l$ with $m$ odd then $\psi(2^kt)=1$ for $k<l$ and
$\psi(2^kt)=0$ for $k\ge l$.  Now the required formula for $E(t)$
follows from the series~(\ref{E-psi}).
\end{proof}

\begin{figure}[ht]
\centering
\mbox{\psfig{file=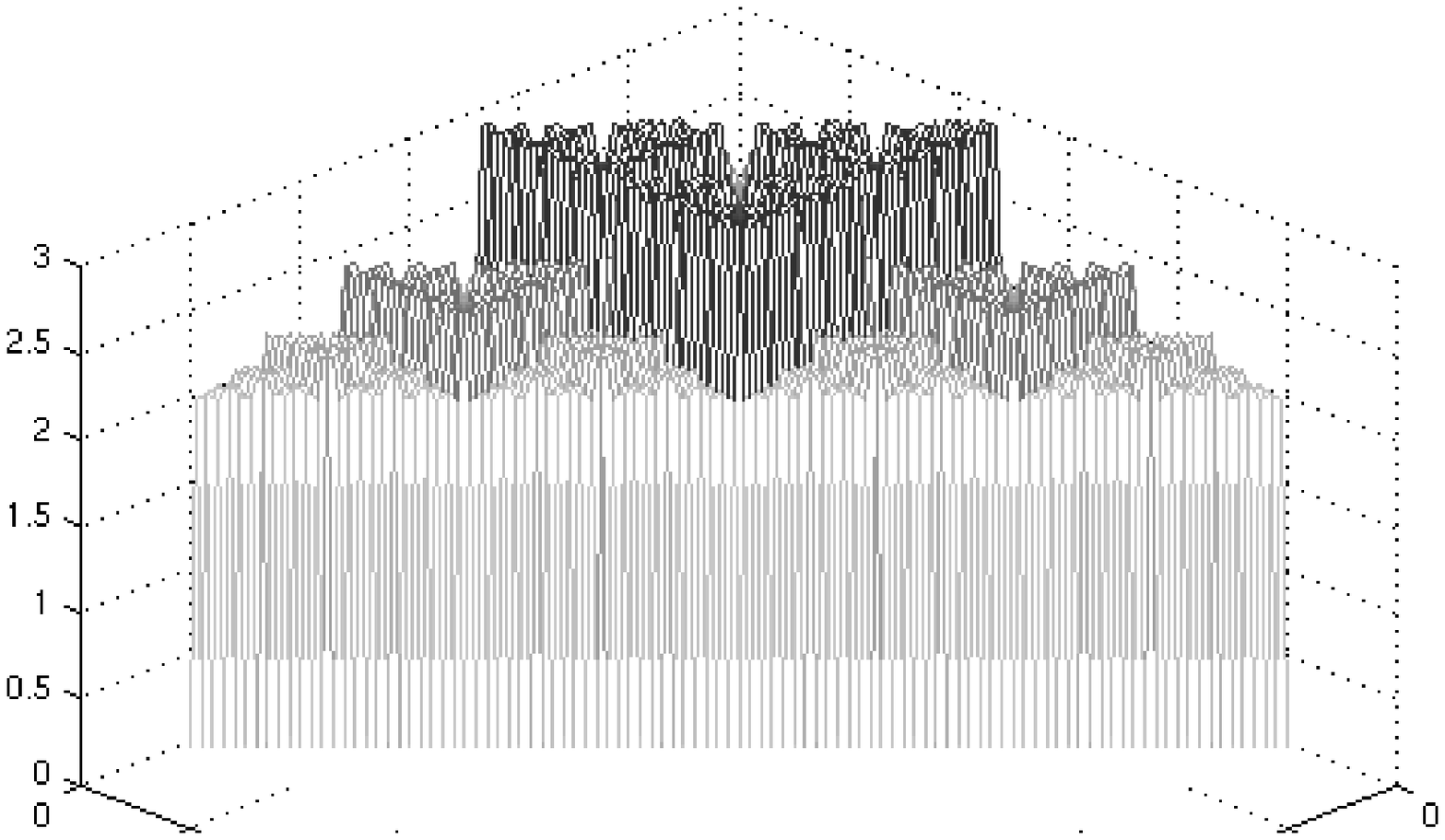,width=5in}}
\caption[]{\mbox{Graph of $z=E(x,y,1-x-y)$ for $0\le x\le 1-y\le 1$} \\
{\tiny  \def\1{\null\hskip.75in}
\null\hskip-.7in
\begin{tabular}{ll}
\1 $\ldots$ tum Tartarus ipse 
		&  Then Tartarus itself goes plunging down\\
 {\em  bis\/} patet in praeceps tantum tenditque sub umbras
		&   In darkness {\em twice\/} as deep as heaven is high\\
   quantus ad aetherium caeli suspectus Olympum
		&     For eyes fixed on etherial Olympus\\
\1	$\cdots$&\1 $\cdots$\\
  Respicit Aeneas subito et sub rupe sinistra
		&  The Heroe, looking on the left, espy'd \\
  moenia lata videt {\em triplici\/} circumdata muro 
		& A lofty Tow'r, and strong on ev'ry side\\
  quae rapidis flammis ambit torrentibus amnis
		& With {\em treble\/} Walls, which Phlegethon surrounds,\\
 Tartareus Phlegethon, torquetque
		&  Whose fiery flood the burning empire bounds:\\
  sonantia saxa 
 &And press'd betwixt the Rocks, the bellowing noise \\
 &  resounds.\\
 & \\
	Vergil, {\it The Aeneid\/} & Translations by Robert Fitzgerald and John Dryden
\end{tabular}}
  }
\label{E2d}
\end{figure}

Unfortunately, in higher dimensions $E$ is not as easy to understand.
A graph of $E$ on the two dimensional simplex is shown in
Figure~\ref{E2d}.

\begin{remark}\label{E-similar}
The graph (Figure~\ref{E2d}) of $E$ suggests that $E$ has some self
similarities.  This is indeed the case as we now briefly indicate.
For each $k\in \{0\cd n\}$ define a map $\theta_k\cn\Delta_n\to
\Delta_n$ by
$$
\theta_k(x):=\frac{e_k+x}{2}
$$
This is the dilation by a factor of $1/2$ centered at $e_k$ and it
maps $\Delta_n$ onto its subset defined by $1/2\le x_k\le
1$. The functional equation~(\ref{H-fcn-eqn}) for $H$ can be
rewritten in the from $H(t/2)=\fp{t/2}+\frac12H(t/2)$.  We leave it as an
exercise for the reader to show this (and $H(t+1/2)=H(t)+1/2$ for
$0< t<1/2$) can be used in the definition of $E$ so show that for any
$x\in \Delta_n$ which is not a vertex that 	
$$
E(\theta_k(x))=1+\frac12 E(x).
$$
Thus if on the space $\Delta_n\times [0,\infty)$ a map $\Theta_k$ is
defined by $ \Theta_k(x,z)=\bigl(\theta_k(x), 1+z/2\bigr)$ then the
graph of $E$ (with the points over the vertices deleted) is invariant
under $\Theta_k$.  Each $\Theta_k$ is the dilation by a factor of
$1/2$ with center $(e_k,2)$.  This explains the self similarities of
the graph of $E$.\qed
\end{remark}

Our next result implies that the upper bound of
Theorem~\ref{thm:convex-bd} is sharp.  Recall that a subset of a
metric space is a $G_\delta$ iff it is a countable intersection of
open sets.

\begin{thm}\label{II3} 
The function $E$ achieves its maximum value of $\kappa(n)$ on an
uncountable $G_\delta$ subset of $\Delta_n$.
\end{thm}

\begin{remark} \label{not-center}
The maximum of $E$ does not occur at the center $(1/(n+1)\cd 1/(n+1))$
of $\Delta_n$.  Given the symmetry of the problem this is a little
surprising.\qed
\end{remark}

\begin{proof}
That $\sup E[\Delta_n]\le \kappa(n)$ follows from~\ref{thm:convex-bd}.
To show the maximum is obtained, let $m=[\log_2(n)]$ so that $n=2^m+r$
with $0\le r<2^m$.  Suppose $x=(x_0\cd x_n)\in \Delta_n$ with each
coordinate~$x_k$ a dyadic irrational.  In particular, if
$x_k=\sum_{j=0}^nx_{kj}/2^j$ then $x_{kj}$ is zero for infinitely
many~$j$ and one for infinitely many~$j$.  Let $M_j(x):=\#\{k: x_{kj}=1\}$.
We claim that $E(x)=\kappa(n)$ provided each coordinate $x_k$ is
a dyadic irrational and
\begin{equation}\label{Mj-def}
M_j(x)=\begin{cases} 0, &j\le m;\\
		n-2r,& j=m+1; \\
		n,   &j\ge m+2. \end{cases}
\end{equation}
Let $K$ be the set of all $x=(x_0\cd x_n)$ that satisfy these two
conditions.  If $x\in K$, then
$$
\sum_{k=0}^nx_k=\sum_{j=0}^\infty\sum_{k=0}^n\frac{x_{kj}}{2^j}=
\frac{M_{m+1}}{2^{m+1}}+n\sum_{j=m+2}^\infty\frac{1}{2^j}
=\frac{n-2r}{2^{m+1}}+\frac{n}{2^{m+1}}=\frac{n-r}{2^m}=1
$$
Thus $x\in \Delta_n$ and so $K\subset \Delta_n$.

To see that $K$ is uncountable (and thus nonempty) let $\la
a_{m+2},a_{m+3},\dots\ra$ be a sequence in $\{0,1\cd n\}$ such that
for every~$k\in \{0,1\cd n\}$, $a_j=k$ for infinitely many~$j$.  We
let $x_{kj}=0$ if $j\le m$ and we let $x_{k\,\, m+1}=1$ for
exactly~$n-r$ many~$k$.  For $j\ge m+2$, let
$$
x_{kj}:=\begin{cases} 1,& a_j\ne j;\\
		0,& a_j=k.\end{cases}
$$
Since each sequence $\la x_{kj}\ra_{j=0}^\infty$ has infinitely many
zeros and ones, each $x_k$ is a dyadic irrational. Thus $x=(x_0\cd
x_n)\in K$.  As there are uncountably many such sequences $\la
a_{m+2},a_{m+3},\dots\ra$ the set $K$ is uncountable.

If $x=(x_0\cd x_n)\in K$ then, using the definition~(\ref{H-def}) of
$H$ and the identity $\sum_{j=m+2}^\infty{j}/{2^j}=(m+3)/2^{m+1}$, we
have
\begin{align*}
E(x)&=\sum_{k=0}^nH(x_k)=\sum_{j=0}^\infty \frac{jM_j}{2^j}
	=\frac{M_{m+1}(m+1)}{2^{m+1}}+n\sum_{j=m+2}^\infty\frac{j}{2^j}\\
&=\frac{(n-2r)(m+1)}{2^{m+1}}+\frac{n(m+3)}{2^{m+1}}
	=\frac{(2n-2r)(m+1)+2n}{2^{m+1}}\\
&=m+1+\frac{n}{2^m}=\kappa(n).
\end{align*}
This shows that $E$ achieves its maximum at all points of $K$.
Finally $\{x\in \Delta_n: E(x)=\kappa(n)\}=\bigcap_{\ell=1}^\infty
E^{-1}\big[(\kappa(n)-1/\ell,\infty)\big]$ and each of the sets
$E^{-1}\big[(\kappa(n)-1/\ell,\infty)\big]$ is open as $E$ is lower
semi-continuous.  Thus $\{x\in \Delta_n: E(x)=\kappa(n)\}$ is
a~$G_\delta$.
\end{proof}

\begin{remark}\label{all-max}
With a little more work it can be shown that $E(x)=\kappa(n)$ if and
only if $x\in K$ with $K$ as above.\qed
\end{remark}

\begin{thm}\label{E-biggest}
The function $E$ is the largest bounded  approximately
convex function on $\Delta_n$ that vanishes on the vertices.  More
precisely, if $h$ is any bounded approximately convex function on
$\Delta_n$ with $h(e_k)\le0$ for $k=0,1\cd n$, then $h\le E$ on
$\Delta_n$.
\end{thm}

 \mnote{\ralph Some rewriting in the corollaries and their proofs.
 Please check}

\begin{cor}\label{Borel-bd}
Let $h\cn \Delta_n\to \R$ be an approximately convex function that is
Borel measurable.  Then for any $x=\sum_{k=0}^nx_ke_k$ the inequality
$$
h(x)\le \kappa(n)+\sum_{k=0}^nx_kh(e_k)
$$
holds. 
In particular, if $h(e_k)\le 0$ for all $k$, then $h\le \kappa(n)$.
\end{cor}

\begin{proof}[Proof of Theorem~\ref{Borel-bd}]
Define $l$ on $\Delta_n$ by $l(x)=\sum_{k=0}^nx_kh(e_k)$.  Then the
function $h(x)-l(x)$ is approximately convex, Borel measurable, and
vanishes on the vertices of $\Delta_n$.  So by replacing $h$ by $h-l$ we
may assume $h$ vanishes on the vertices of $\Delta_n$  it will be
enough to show $h\le \kappa(n)$ on $\Delta_n$.  We do this by
induction on $n$.  For $n=1$ it follows from results of Ng and
Nikodem~\cite[Cor.~1 and Thm~2]{Ng-Nikodem} that $h$ is bounded above.
But then $h\le \kappa(1)=2$ by Theorem~\ref{thm:convex-bd}.
Now let $n\ge 2$ and assume the result
holds for all simplices with dimension $<n$.  Consider $\Delta_{n-1}$
as a face of $\Delta_n$ in the natural way
($\Delta_{n-1}=\conv\{e_0\cd e_{n-1}\}\subset \conv\{ e_0\cd e_n\}$).
Then by the induction hypothesis $h\big|_{\Delta_{n-1}}\le
\kappa(n-1)$.  Now any point $x\in \Delta_n$
has a representation as $x=(1-t)e_n+ty$ where $y\in \Delta_{n-1}$ and
$t\in [0,1]$.  But then the one dimensional result
(applied to the restriction of $h$ to the segment between $e_0$ and
$y$ where we note that this restriction is Borel and thus Lebesgue
measurable) implies
\begin{align*}
h(x)&= h((1-t)e_n+ty) \le 2+(1-2)h(e_n) + th(y)\\
& \le 2+0+ t\kappa(n-1)  \le 2+\kappa(n-1). 
\end{align*}
Thus $h$ is bounded above on $\Delta_n$.  But then we can use
Theorem~\ref{thm:convex-bd} and reduce the bound to $\kappa(n)$.  This
completes the proof.
\end{proof}

\begin{cor}\label{cor:E-combo}
Let $U\subseteq \R^n$ be a convex set and let $h\cn U\to \R$ be either
Borel measurable or bounded above on compact subsets of $U$.  Then for
any $m\le n$, points $x_0\cd x_m\in U$ and $(\alpha_0\cd \alpha_m)\in
\Delta_m$, we have
\begin{align*}
h(\alpha_0x_0+\cdots+ \alpha_mx_m)&\le E(\alpha_0\cd \alpha_m)+
\alpha_0h(x_0) +\cdots +\alpha_mh(x_m)\\
&\le \kappa(m)+ \alpha_0h(x_0) +\cdots +\alpha_mh(x_m).
\end{align*}
\end{cor}

\begin{proof}
Define $f\cn \Delta_m\to \R$ by $f(\alpha_0\cd \alpha_m):=
h(\alpha_0x_0+\cdots+ \alpha_mx_m)-\big(\alpha_0h(x_0) +\cdots
+\alpha_mh(x_m)\big)$.  Then $f$ is approximately convex and bounded
above on $\Delta_m$ or is Borel measurable on $\Delta_m$.  As $f$
vanishes on the vertices of $\Delta_m$ either Theorem~\ref{E-biggest}
or Corollary~\ref{Borel-bd} implies $f\le E\le
\kappa(m)$ on $\Delta_m$.  This is equivalent to the conclusion of the
corollary.
\end{proof}

We start the proof of Theorem~\ref{E-biggest} by extending the idea of
the dyadic support from $\dy$ to $\dy_n$.  If
$x=\sum_{k=0}^n\(\sum_{j=0}^Nx(j,k)/2^j\)e_k\in \dy_n$ (here
$x(j,k)\in \{0,1\}$) then set
\begin{equation}\label{dy-supp}
\supp x:=\{(j,k): x(j,k)=1\}.
\end{equation}
The following is trivial to prove using Lemma~\ref{H-exact} and the
definition of $E$ in terms of $H$.

\begin{lemma}\label{E-exact}
If $x,y\in \dy_n$ and $(\supp x)\cap (\supp y)=\nothing$ then
$$
E\Bigl(\frac{x+y}{2}\Bigr)=\frac{E(x)+E(y)}{2}+1.
$$
\vskip-15pt\qed
\end{lemma}

\begin{lemma}\label{supp-split}
If $x\in \dy_n$ and $x\notin\{e_0\cd e_n\}$, then there are $y,z\in
\dy_n$ so that $x=(y+z)/2$ and $(\supp y)\cap (\supp z)=\nothing$.
\end{lemma}

\begin{proof}
Letting  $x=\sum_{k=0}^n\(\sum_{j=0}^Nx(j,k)/2^k\)e_k$
It suffices to show that there are nonempty sets $A$, $B$ so 
that $A\cap B=\nothing$ and
$$
\sum_{(j,k)\in A}\frac{x(k,j)}{2^j}
	 =\frac12=\sum_{(j,k)\in B}\frac{x(k,j)}{2^j}.
$$
For then if $a=\sum_{(j,k)\in A}{x(k,j)}/{2^{j-1}}e_k$ and
$b=\sum_{(j,k)\in B}{x(k,j)}/{2^{j-1}}e_k$ we have $a,b\in\dy_n$,
$(\supp a)\cap (\supp b)=\nothing$ and $x=(a+b)/2$.

We first prove by induction on $\sum_{j=1}^Na_j$ that if $a_1\cd a_N$ are
positive integers so that $\sum_{j=1}^Na_j/2^j=1$ then
there are $b_j,c_j\in\Nat$ such that
$\sum_{j=1}^Nb_j/2^j=\sum_{j=1}^Nc_j/2^j=1/2$.  Note that $a_N$ is
even (otherwise $2^{-N}\sum_{j=1}^N2^{N-j}a_j$ would not sum
to~$1$) and so $a_N-2\ge 0$.  Therefore
$$
\sum_{j=1}^{N-2}\frac{a_j}{2^j}+\frac{a_{N-1}+1}{2^{N-1}}+\frac{a_N-2}{2^N}=1.
$$
Since $\sum_{j=1}^{N-2}+(a_{N-1}+1)+(a_N-2)=\sum_{j=1}^Na_j-1$ we may
apply the induction hypothesis, which yields the claim.

Now let $x\in \dy_n$ be as above.  Let $a_j:=\#\{k:x(j,k)=1\}$.  Then
$\sum_{j=1}^Na_j/2^j=1$.  Therefore we have $a_j=b_j+c_j$ as above.
Then splitting each of the sets $\{k:x(j,k)=1\}$ into two disjoint
sets $A_j$ and $B_j$ with $\#(A_j)=b_j$ and $\#(B_j)=c_j$ we let
$A:=\cup_{j=1}^N A_j$ and $B:=\cup_{j=1}^NB_j$.  This completes the
proof.
\end{proof}

\begin{prop}\label{dy-E-biggest}
Let $h$ be any approximately convex function on $\Delta_n$ (not necessarily
bounded above) such that $h(e_k)\le 0$ for $0\le k\le
n$.  Then $h(x)\le E(x)$ for all $x\in \dy_n$.
\end{prop}

\begin{proof}
The proof is by induction on $m=\#(\supp x)$.  If $m=1$ then  $x=e_k$
for some $k$ and $h(e_k)\le 0 =E(e_k)$.  Now assume that $h(x)\le
E(x)$ for all $x$ with $\#(\supp x)\le m-1$ and let $\supp x=m$.  By
Lemma~\ref{supp-split} we can write $x=(y+z)/2$ with $\#(\supp y),
\#(\supp z)\le m-1$.   Using the induction hypothesis and
Lemma~\ref{E-exact} 
$$
h(x)=h\Bigl(\frac{y+z}{2}\Bigr)\le \frac{h(y)+h(z)}{2}+1
\le \frac{E(y)+E(z)}{2}+1=E(x).
$$
\vskip-15pt
\end{proof}

The following lets us pass from knowing inequalities for $E$ on
$\dy_n$ to proving them on $\Delta_n$.

\begin{lemma}\label{E-lim-dy}
If $x\in \Delta_n$ then there is a sequence $\la x(r)\ra_{r=1}^\infty$
from $\dy_n$ so that $\lim_{r\to\infty}x(r)=x$ and
$\lim_{r\to\infty}E(x(r))=E(x)$. 
\end{lemma}

\begin{proof} 
Write $x=\sum_{k=0}^nx_ke_k$.  By reordering we can assume for some
$\ell\in\{0\cd n\}$ that $x_k\in \dy$ for $0\le k\le \ell$ and
$x_k\notin \dy$ for $\ell+1\le k\le n$.  For $0\le k\le \ell$ set
$x_k(r)=x_k$ for all $r$.  As $\sum_{k=0}^nx_k=1$ and $\sum_{k=0}^\ell
x_k\in \dy$ (as $x_k\in \dy$ for each $x_k$ in this sum) the sum
$\delta:=\sum_{k=\ell+1}^n x_k=1-\sum_{k=0}^\ell x_k $ will also be a
dyadic rational.  Let
$\Delta_{n-\ell-1}(\delta)=\{\sum_{k=\ell+1}^n\alpha_ke_k:\alpha_k\ge0,
\sum_{k=\ell+1}^n\alpha_k=\delta\}$ and
$\dy_{n-\ell-1}(\delta)=\{\sum_{k=\ell+1}^n\alpha_ke_k:\alpha_k\in\dy,
\sum_{k=\ell+1}^n\alpha_k=\delta\}$.  Then $\dy_{n-\ell-1}(\delta)$
will be dense in $\Delta_{n-\ell-1}(\delta)$ so there is a sequence
$y(r)=\sum_{k=\ell+1}^n y_k(r)e_k$ with $\lim_{r\to\infty}y(r)=y$.
Set $x_k(r)=y_k(r)$ for $\ell+1\le k\le n$.  Then $x_k(r)=x_k\in \dy$
for $0\le k\le \ell$ and $\lim_{r\to \infty}x_k(r)=x_k\notin \dy$ for
$\ell+1\le k\le n$. Set $x(r)=\sum_{k=0}^nx_k(r)e_k$.  Then
$x(r)\in\dy_n$ and $\lim_{r\to\infty}x(r)=x$.  We now use the
definition of $E$ in terms of $H$ and the fact that $H$ is continuous at all
dyadic irrationals (Proposition~\ref{H-series}) to obtain
$$
\lim_{r\to\infty}E(x(r))=\sum_{k=0}^\ell
H(x_k)+\lim_{r\to\infty}\sum_{k=\ell+1}^nH(x_k(r))
=\sum_{k=0}^n H(x_k)=E(x).
$$
\vskip-15pt
\end{proof}

\begin{proof}[Proof of Theorem~\ref{E-biggest}]  Let $E_{\Delta_n,0}$
extremal approximately convex function on $\Delta_n$ that takes the
values~$0$ on the vertices (cf.~(\ref{E-K-def})).  We wish to
show $E=E_{\Delta_n,0}$.  The inequality $E\le E_{\Delta_n,0}$ follows from the
definition of $E_{\Delta_n,0}$, so it is enough to prove $E_{\Delta_n,0}\le E$.  By
Lemma~\ref{E-lim-dy} there is a sequence $x(r)\in\dy_n$ such that
$\lim_{r\to\infty}x(r)=x$ and $\lim_{r\to\infty}E(x(r))=E(x)$.  
By Lemma~\ref{dy-E-biggest} $E_{\Delta_n,0}(x(r))\le E(x(r))$.  By
Theorem~\ref{extreme-lsc} the function $E_{\Delta_n,0}$ is lower
semi-continuous.  Therefore
$$
E_{\Delta_n,0}(x)\le \liminf_{r\to \infty}E_{\Delta_n,0}(x(r))
	\le \lim_{r\to\infty}E(x(r))=E(x).
$$
\vskip-15pt
\end{proof}

\subsection{Extremal approximately convex functions on convex polytopes}
\label{sec:polytopes}

Let $K\subset \R^n$ be a compact convex set with extreme points $V$
and let $\phi\cn V\to \R$ be bounded.  In \S\ref{sec:lsc} we defined
the extremal approximately convex function $E_{K,\phi}$ with extreme
values $\phi$ but without being explicit about how to compute it.
In~\S\ref{sec:E} we gave a very explicit description of
$E=E_{\Delta_n,0}$, the extremal approximately convex function on the
simplex.  Here we show that when $K$ is a polytope (that is the convex
hull of a finite number of points) then $E_{K,\phi}$ can be expressed
directly in terms of $E_{\Delta_m,0}$ for some $m$.  We first
establish some elementary properties of approximately convex functions
under affine maps.

\begin{prop}\label{map-props}
Let $A\subset \R^m$ and $B\subset \R^n$ be convex sets and $T\cn
\R^m\to \R^n$ an affine map.  
\begin{enumerate}
\item If $T[A]\subseteq B$ and $f$ is an approximately convex function
on $B$ then $T^*f(x):=f(T(x))$ is an approximately convex function on
$A$. 
\item If $T[A]\supseteq B$ and $h$ is an approximately convex function
on $A$ which is bounded from below 
then $T_*h(y):=\inf_{T(x)=y}h(x)$ is approximately convex on $B$. 
\item Both $T^*$ and $T_*$ are order preserving.  That is $f_1\le f_2$
and $h_1\le h_2$ pointwise implies $T^*f_1\le T^*f_2$ and $T_*h_1\le
T_*h_2$ pointwise.
\item If $T[A]=B$, $h$ is approximately convex and bounded below on
$A$ and $f$ is approximately convex and bounded below on $B$, then
$T^*T_*h\le h$ and $T_*T^*f=f$.
\end{enumerate}
\end{prop}

\begin{proof}
This is just a chase through the definitions of $T^*$ and $T_*$.
\end{proof}

Let $K$ be a convex polytope in $\R^n$ with extreme points $V=\{v_0\cd
v_m\}$ and extreme values given by $\phi\cn V\to \R$.  and let $a_0\cd
a_m$ be real numbers.  We wish to find the largest approximately
convex function $F$ on $K$ so that $F(v_k)=\phi(v_k)$ for $0\le k\le
m$.  Toward this end let $E=E_{\Delta_m,0}$ be the extremal
approximately convex function on the simplex~$\Delta_m$ and define
$E_{\Delta_m,\phi}$ on $\Delta_M$ by
$$
E_{\Delta_m,\phi}(x)=E_{\Delta_m,\phi}(x_0\cd
x_m):=E_{\Delta_m,0}(x_0\cd x_m)+\sum_{k=0}^mx_k\phi(a_k).
$$
(This is a slight misuse of notation as $\phi$ is a function on the
extreme points $V$ of $K$ rather than the set of extreme points
$\{e_0\cd e_m\}$ of $\Delta_m$.)  Then, as $x\mapsto
\sum_{k=0}^mx_ka_k$ is affine, the function $E_{\Delta_m,\phi}$ is
approximately convex on $\Delta_m$ and satisfies
$E_{\Delta_m,\phi}(e_k)=\phi(v_k)$.  Moreover $E_{\Delta_m,\phi}$ is
the extremal approximately convex function on $\Delta_m$ taking on
these values on the vertices in the sense that if $f\cn \Delta_m\to
\R$ is approximately convex and bounded above, lower semi-continuous,
and $f(e_k)\le \phi(v_k)$ then $f(x)\le E_{\Delta_m,\phi}(x)$ for all $x\in
\Delta_m$.

Returning to our extremal problem there is a unique affine map $T\cn
\Delta_m\to K$ such that $T(e_k)=v_k$ for $0\le k \le m$.  Then
$T[\Delta_m]=K$.  Define $F_{K,}\cn K\to \R$ by
$$
F_{K,\phi}:=T_*E_{\Delta_m,\phi}.
$$
Then another definition chase shows $F_{K,\phi}(v_k)=a_k$.  

\begin{thm}\label{gen-extremal}
Using the notation above, the extremal approximately continuous
function on the polytope $K$ with extreme values $\phi$ is
$$
E_{K,\phi}:=T_*E_{\Delta_0,\phi}.
$$
The function $E_{K,\phi}$ is lower semi-continuous.
\end{thm}

\begin{proof}
Let $f\cn K\to \R$ be approximately convex, bounded, and and satisfy
$f(v_k)\le \phi(v_k)$.  Then the function $T^*f$ on $\Delta_m$ is
approximately convex, bounded, and $T^*f(e_k)=f(v_k)\le \phi(v_k)$.
Therefore $T^*f\le E_{\Delta_m,\phi}$.  But then $f=T_*T^*f\le
T_*E_{\Delta_m,\phi}$ which proves $T_*E_{\Delta_m,\phi}=E_{K,\phi}$.
The lower semi-continuity of $E_{K,\phi}$ follows from
Theorem~\ref{extreme-lsc}.
\end{proof}

\subsection{A stability theorem of Hyers-Ulam type}
\label{sec:Hyers-Ulam}
Here we give a stability result for approximately convex functions
related to and motivated by a theorem of Hyers and
Ulam~\cite{Hyers-Ulam:convex}.  The idea is that an approximately
convex function is close (in the uniform norm) to some convex
function.

\begin{thm}\label{thm:stable}
Assume that $U\subseteqq \R^n$ is convex, $\e>0$, and that $f\cn U\to
\R$ is bounded above on compact sets and satisfies
\begin{equation}\label{eps-convex}
f\Bigl(\frac{x+y}{2}\Bigr)\le \frac{f(x)+f(y)}{2}+\e.
\end{equation}
Then there exist convex functions $g, g_0\cn U\to \R$ such that
\begin{equation}\label{ineq-stable}
f(x) \le g(x)\le f(x)+\kappa(n)\e\quad \text{and}\quad 
|f(x)-g_0(x)|\le \frac{\kappa(n)}{2}\e
\end{equation}
for all $x\in U$.  The constant $\kappa(n)$ is the best possible
constant in these inequalities.
\end{thm}

\begin{proof}
By replacing $f$ by $\e^{-1}f$ we may assume $\e=1$ so that $f$ is
approximately convex.
Following
Hyers and Ulam~\cite[p.~823]{Hyers-Ulam:convex} or
Cholewa~\cite[pp.~81--82]{Cholewa:stability} set $W:=\{(x,y)\in
\R^n\times \R: y\ge f(x)\}$ and define $g$ by
$$
g(x):=\inf\{y: (x,y)\in \conv(W)\}.
$$
We now show that $g$ does not take on the value $-\infty$.  If
$(x,y)\in \conv(W)$ then by  Carath\'eodory's Theorem there exist $n+2$ points
$(x_0,y_0)\cd (x_{n+1},y_{n+1})\in W$ and $(\alpha_0\cd \alpha_n)\in
\Delta_{n+1}$ such that $(x,y)=\sum_{k=0}^{n+1}\alpha_k(x_k,y_k)$.
Therefore by Corollary~\ref{cor:E-combo}
\begin{align*}
f(x)&=f\biggl(\sum_{k=0}^{n+1}\alpha_kx_k\biggr)\le \kappa(n+1)+
\sum_{k=0}^{n+1}\alpha_kf(x_k)\\
& \le \kappa(n+1)+\sum_{k=0}^{n+1}\alpha_ky_k=\kappa(n+1)+y.
\end{align*}
Thus $y\ge f(x)-\kappa(n+1)$ which implies $g(x)\ge
f(x)-\kappa(n+1)>-\infty$.  

>From the definition it is clear that $g(x)\le f(x)$ and that $g(x)$ is
convex.  To see that $f(x)\le g(x)+\kappa(n)$ let $\delta>0$ and
choose $y$ so that $(x,y)\in \conv(W)$ and $y<g(x)+\delta$.  Then as
above there are $n+2$~points $(x_0,y_0)\cd (x_{n+1},y_{n+1})\in W$
with $(x,y)\in W$ and such that $(x,y)\in \Delta:=\conv(\{(x_0,y_0)\cd
(x_{n+1},y_{n+1})\})$.  Let $\ol{y}:=\min\{\eta: (x,\eta)\in
\Delta\}$.  Then $(x,\ol{y})$ is on the
boundary of $\Delta$ and so it is a convex combination of $n+1$ of the
points $(x_0,y_0)\cd (x_{n+1},y_{n+1})$, say
$(x,\ol{y})=\sum_{k=0}^n\alpha_k(x_k,y_k)$ with $(\alpha_0\cd
\alpha_n)\in \Delta_n$.  Then a calculation like one showing that
$g(x)>-\infty$ (but with $n+1$ replacing $n+2$) yields that $f(x)\le
\ol{y}+\kappa(n)\le g(x)+\delta+\kappa(n)$.  As $\delta>0$ was
arbitrary this implies $f(x)\le g(x)+\kappa(n)$.

Letting $g_0(x)=g(x)+\kappa(n)/2$ we have $|f(x)-g_0(x)|\le
\kappa(n)/2$. 

Finally to see that the constants in question are sharp consider the
almost convex function $E\cn \Delta_n\to \R$ which has $\max
E=\kappa(n)$.  Then the largest convex function $g$ on $\Delta_n$
with $g\le E$ is $g(x)\equiv0$.  Likewise $g_0(x)\equiv\kappa(n)/2$
has $|E(x)-g_0(x)|\le \kappa(n)/2$ and no other convex function on
$\Delta_n$ gives a better estimate.
\end{proof}

\subsection{Examples of approximately convex functions}\label{sec:examples}
Here we give examples showing that the hypothesis of our results are
necessary.
\begin{example}\label{sub-convex}
Let $f(t)$ be any approximately sub-affine function on $[0,1]$.  Then (as in
the proof of Proposition~\ref{E-basics}) the function
$F(x):=f(x_0)+f(x_1)+\cdots +f(x_n)$ defined on the simplex $\Delta_n$ 
will be approximately convex.  Using the
function $f(t)=t\log_2(1/t)$ shows that for example 
$F(x):=\sum_{k=0}^nx_k\log_2(1/x_k)$ is approximately convex
(cf.~Lemma~\ref{x*log(x)-ineq}). As a slight
generalization of this if $f_0\cd f_n$ are all approximately sub-affine then 
$F_1(x)=f_0(x_0)+f_1(x_1)+\cdots +f_n(x_n)$ is approximately convex.\qed
\end{example}

\begin{example}\label{ex:bounded}
Let $C$ be any convex subset of any normed vector space and let
$\phi\cn C\to [0,1]$.  Then $\phi((x+y)/2)\le 1\le
(\phi(x)+\phi(y))/2+1$ so $\phi$ is approximately convex.  There is no
assumption on $\phi$ other than the bounds $0\le \phi\le 1$.  Thus
$\phi$ need not be continuous or measurable.  So approximate convexity
by itself does not imply any type of regularity of the function.\qed
\end{example}

\begin{example}\label{ex:Q-linear} 
View $\R^{n+1}$ as a vector space over the rational numbers $\Q$ and
let $\mathcal B$ be a Hamel basis for $\R^{n+1}$ over $\Q$.  Let $h\cn
\R^{n+1}\to \R$ obtained by first mapping $\mathcal B$ to $\R$ and then
extending to $\R^{n+1}$ by linearity.  We can choose $\mathcal
B\subset \Delta_n=\conv\{e_0\cd e_n\}$ (with $e_0\cd e_n$ the standard
basis of $\R^{n+1}$) and $h$ so that $h[\mathcal B]$ is dense in
$\R$.  Therefore $h$ is unbounded on $\Delta_n$. 

To get an example more closely related to Theorem~\ref{thm:convex-bd}
let $h$ be as just defined but chosen in such a way that $h(e_i)=0$
for $0\le i\le n$ and set $ h_0(x):=\max \{h(x),0\}$.  Then for
$A:=\{e_0\cd e_n\}$ we have $\Delta_n=\conv(A)$, $h_0$ is bounded from
below, and $h_0\equiv 0$ on $A$.  But $h_0$ is not bounded from above
on $\Delta$.  This shows the assumption that $h$ be bounded from above
in Theorem~\ref{thm:convex-bd} is necessary.  A similar example appears
in the paper of Cholewa~\cite[\S3]{Cholewa:stability}.  \qed
\end{example}

\begin{example}\label{ex:not-meas}
As an extension of the last example let $\Delta^{k}_{n-1}$ for $0\le
k\le n$ be the ($(n-1)$-dimensional) faces of $\Delta_n$.  For each
$k$ choose an unbounded approximately convex function
$h_k:\Delta^{k}_{n-1}\to [0,\infty)$ that vanishes on the vertices of
$\Delta_{n-1}^k$ (possible by the last example).  Let $h\cn\Delta_n\to
[0,\infty)$ be $h(x)=0$ on the interior of $\Delta_n$ and for each face
$h\big|_{\Delta_{n-1}^k}=h_k$. (A little care must be taken  in the
choice of the $h_k$'s to ensure that these restrictions agree on the
intersections of the faces. This is not hard to arrange and we leave
the details to the reader.) Then as the boundary of $\Delta_n$
(which is $\bigcup_{k=0}^n\Delta_{n-1}^k$) is a set of measure zero
the function $h$ is Lebesgue measurable on $\Delta_n$, but is not
Borel measurable.  This shows that the hypothesis of
Theorem~\ref{Borel-bd} can not be weakened from Borel measurable to
Lebesgue measurable.~\qed
\end{example}

\section{The Size of the Convex Hull of an Approximately Convex
Set}\label{sec:sets} 

In this section we apply our results on approximately convex functions to the
problem of giving {\it a priori\/} bounds on the size of convex hull
of an approximately convex set.

\subsection{General upper bounds}
\label{sec:gen-upper}

We now apply our results to the geometric problem of computing the
size of the convex hull.

\begin{thm}\label{gen-upper}
Let $\|\cdot\|$ be any norm on $\R^n$ and let $A\subset \R^n$ be a set
that is approximately convex in this norm.  Let $b\in \conv(A)$ so that for
some $a_0\cd a_m \in A$ with $m\le n$ we have
$b=\sum_{k=0}^m\alpha_ka_k$ where $(\alpha_0\cd \alpha_n)\in
\Delta_m$, then
\begin{equation}\label{dist-bd}
\dist(b,A)\le E(\alpha_0\cd \alpha_m)\le \kappa(m)\le \kappa(n).
\end{equation}
(In the terminology of Theorem~\ref{midpoint} this implies that 
$C_{\|\cdot\|}\le
\kappa(n)$.) 
\end{thm}

\begin{remark}\label{rescale}
For bounded sets this result can be restated in a dilation invariant
fashion that does not involve approximately convex sets in its
statement: If $A\subset \R^n$ is bounded set and $b\in \conv(A)$ so
that $b=\sum_{k=0}^m\alpha_ka_k$ as in the statement of the theorem,
then
$$
\dist(b,A)\le  E(\alpha_0\cd \alpha_m)d_H(\frac12(A+A),A)
\le \kappa(m)d_H(\frac12(A+A),A).
$$
The results below have similar dilation invariant versions.\qed
\end{remark}

\begin{proof}
Define a function $f\cn \Delta_m\to [0,\infty)$ by 
$$
f(\beta_0\cd \beta_m):= \dist\bigg(\sum_{k=0}^m\beta_ma_k,A\bigg).
$$
Then as the function function $x\mapsto \dist(x,A)$ on $\R^n$ is an
approximately convex function and the map $(\beta_0\cd \beta_m)\mapsto
\sum_{k=0}^n\beta_ka_k$ is affine the function $f$ is approximately convex
and it is clearly continuous.  Also $f$ vanishes on the vertices of
$\Delta_m$.  Therefore by Theorem~\ref{E-biggest} the bound
$f(\beta_0\cd \beta_n)\le E(\beta_0\cd \beta_n)$ holds.  But this
implies~(\ref{dist-bd}).
\end{proof}

Recall that a subset $A\subset \R^n$ is {\bi convexly connected\/} iff
there is no hyperplane $H$ of $\R^n$ so that $A$ meets both half
spaces determined by $H$ but does not meet $H$.  Each subset $A$
decomposes uniquely into convexly connected components.

\begin{thm}\label{connect-bd}
Let $\|\cdot\|$ be a norm on $\R^n$ and let $A\subset \R^n$ which is
approximately convex in this norm.  Assume that either $A$ has at most $n$
connected components or $A$ is compact and has at most $n$ convexly
connected components.  Then any $b\in \conv(A)$ satisfies
$\dist(b,A)\le \kappa(n-1)$.
\end{thm}

\begin{proof}
In either of the two cases there is a refinement of Carath\'eodory's
Theorem (cf.~\cite{Hanner-Radstrom}) which implies that $b$ is a convex
combination of $n$~points $a_0\cd a_{n-1}$ points of $A$.  Then
Theorem~\ref{gen-upper} with $m=n-1$ implies $\dist(b,A)\le
\kappa(n-1)$.
\end{proof}

In a normed space we will use the notation  $B_R(x_0)$ for the closed ball
of radius $R$ about $x_0$.

\begin{prop}\label{trap0}
Let $\|\cdot\|$ be a norm on $\R^n$ and $A\subset \R^n$ a closed
subset of $\R^n$.  Assume that $x_0\in\R^n\setminus A$ is a point
where the function $x\mapsto \dist(x,A)$ has a local maximum.  Set
$R:=\dist(x_0,A)$ and let $A_1:=B_R(x_0)\cap A$ be the points of
$A$ at a distance $R$ from $x_0$.  Then there are points $a_0\cd
a_k\in A_1$ with $k\le n$ and norm one linear functionals $\lambda_0\cd
\lambda_k\in \R^{n*}$ so that $\lambda_i(a_i-x_0)=R$ (i.e. $\lambda_i$ norms
$a_i-x_0$) and with $0\in \conv\{\lambda_0\cd \lambda_k\}$.  
\end{prop}

\begin{proof} 
By translation and rescaling we may assume $x_0=0$ and $R=1$.
Let $S:=\{u\in \R^n:\|u\|=1\}$ be the unit sphere of the norm
$\|\cdot\|$.  Let $\|\cdot\|^*$ be the dual norm on $\R^{n*}$ and
$S^*$ the unit sphere of $\|\cdot\|^*$.  For any subset $C\subset
\R^n$ let $N^*(C)$ be the set of linear functionals that norm some
member of $C$.  Explicitly 
$
N^*(C):=\{\lambda\in S^*: \lambda(c)=\|c\|\text{ for some } c\in C\}.
$
If $C$ is compact then $N^*(C)$ is also compact. (For if $\la
\lambda_\ell\ra_{\ell=1}^\infty$ is a sequence from $N^*(C)$ then (as
$S^*$ is compact) by going to a subsequence we can assume that
$\lambda_\ell\to \lambda$ for some $\lambda\in S^*$.  For each $\ell$
there is a $c_\ell \in C$ with $\lambda_\ell(c_\ell)=\|c_\ell\|$.  By
compactness of $C$ and again going to a subsequence we  assume
$c_\ell\to c$ for some $c\in C$.  But then
$\lambda(c)=\lim_{\ell\to0}\lambda_\ell(c_\ell)
=\lim_{\ell\to0}\|c_\ell\|=\|c\|$
which shows $\lambda\in N^*(C)$.  Thus any sequence from $N^*(C)$
contains a subsequence that converges to a point of $N^*(C)$ and
therefore $N^*(C)$ is is compact.)

Let $d_H(\cdot,\cdot)$ be the Hausdorff distance defined on the
compact subsets of $\R^n$.  View the map $C\mapsto N^*(C)$ as a map
from the set of compact subsets of $\R^n$ to the set of compact
subsets of $S^*$.  Then we claim this map is {\bi sub-continuous} in
the sense that if $d_H(C_\ell, C)\to 0$ and $K\subseteq S^*$ is a
cluster point of the sequence $\la N^*(C_\ell)\ra_{\ell=1}^\infty$
then $K\subseteq N^*(C)$.  To see this note as $K$ is a cluster point
of $\la N^*(C_\ell)\ra_{\ell=1}^\infty$ by going to a subsequence we
can assume $N^*(C_\ell)\to K$.  Choose $\lambda\in K$.  Then we can
choose $\lambda_\ell\in N^*(C_\ell)$ in such a way that
$\lambda_\ell\to\lambda$.  From the definition of $N^*(C_\ell)$ there
is a $c_\ell\in C_\ell$ so that $\lambda_\ell(e_\ell)=\|c_\ell\|$.  By
yet again going to a subsequence it can be assumed $c_\ell\to c$ for
some $c\in C$.  But then a calculation like the one showing $N^*(C)$
is compact yields $\lambda(c)=\|c\|$.  Thus $\lambda\in N^*(C)$.  As
$\lambda$ was any element of $K$ this shows $K\subset N^*(C)$ as
claimed.

Returning to the proof of Proposition~\ref{trap0}.  For $r\ge 1$ let
$A_r:=\{a\in A:\|a\|\le r\}$.  Then, as in the statement of the
proposition, $A_1$ is the set of points of $A$ at a distance
exactly~$1$ from $0$ and so the conclusion of the proposition is
equivalent to $0\in \conv(N^*(A_1))$ (for if $0$ is a convex
combination of elements of $N^*(A_1)$ then the number of elements can
be reduced to $n+1$ by Carath\'eodory's Theorem).  Assume, toward a
contradiction, that $0\notin \conv(N^*(A_1))$.  Then $N^*(A_1)$ is
compact and thus $\conv(N^*(A_1))$ is also compact.  Therefore the
distance from $\conv(N^*(A_1))$ to $0$ is positive, say $2\delta$.  As
$1\le r\le s$ implies $A_1\subseteq A_r\subseteq A_s$ and
$\bigcap_{r\ge 1}A_r=A_1$ it is not hard to see that $\lim_{r\searrow
1}d_H(A_r,A_1) =0$.  Thus by the sub-continuity of $N^*$ there is an
$r_0>1$
so that the set $N^*(A_{r_0})$ has Hausdorff distance $<\delta$ from
some subset $K$ of $N^*(A_1)$.  This implies the Hausdorff distance
between $\conv(N^*(A_{r_0}))$ and $\conv(K)$ is $<\delta$ and as
$K\subset N^*(A_1)$ this implies $\dist(0,N^*(A_{r_0}))\ge \delta$.
Thus there is a a linear functional on $\R^{n*}$ that separates
$N^*(A_{r_0})$ from $0$.  As the linear functionals on $R^{n*}$ are
the point evaluations there is a unit vector $u_0\in S$ and $\e>0$ so
that for all $\lambda\in N^*(A_{r_0})$ the inequality $\lambda(u_0)\le
-\e$ holds.  Therefore for any $b\in A_{r_0}$ we have a $\lambda\in
N^*(A_{r_0})$ that norms $b$ and so for all $t>0$
$$
\|b-tu_0\|\ge \lambda(b-tu_0)=\|b\|-t\lambda (u_0) \ge 1+\e t
$$
and so $\dist(tu_0,A_{r_0})\ge 1+\e t$ for all $t\ge0$. 
Suppose that $\|x\| < (r_0-1)/2$.  Then $\dist(x,A)\le
\dist(0,A)+\|x\|< 1+ (r_0-1)/2= (r_0+1)/2$. Suppose that $a\in A$ and
that $\|a\|>r_0$.  Then $\|a-x\|>r_0-\|x\| > (1+r_0)/2 > \dist(x,A)$.
Thus, $\dist(x,A)=\dist(x,A_{r_0})$.
In particular this implies that for
$0<t<(r_0-1)/2$ that $\dist(tu_0,A)=\dist(tu_0,A_{r_0})\ge 1+\e t>1$.
This contradicts that $\dist(\cdot,A)$ has a local maximum at $x=0$
and completes the proof.
\end{proof}

\subsection{General lower bounds.}
\label{sec:sharp}

The following result shows that the estimate of
Theorem~\ref{gen-upper} is sharp for all $m\le n-1$ and that
Theorem~\ref{connect-bd}, Theorem~\ref{thm:Euc-bd} and
Theorem~\ref{2D-bds} are all sharp.

\begin{thm}\label{II4}  
Let $\|\cdot\|$ be any norm on $\R^n$ with $n\ge 2$ and let
$\alpha=(\alpha_0\cd \alpha_{n-1})\in \Delta_{n-1}$.  Then, for any
$\e>0$, there is a compact connected approximately convex set $A\subset \R^n$
and a point $b\in \conv(A)$ so that $b=\sum_{k=0}^{n-1}\alpha_ka_k$,
with $a_k\in A$, so that $\dist(b,A)\ge E(\alpha_0\cd
\alpha_{n-1})-\e$.  In particular, since $\sup_{x\in
\Delta_{n-1}}E(x)=\kappa(n-1)$ (cf.~\ref{II3}),
for the proper choice of $\alpha$ it
follows that there is a compact connected approximately convex set $A\subset \R^n$
and a point $b\in A$ so that $\dist(b,A)\ge \kappa(n-1)$. (In the
terminology of Theorem~\ref{midpoint} this implies that $C_{\|\cdot\|}\ge
\kappa(n-1)$.)
\end{thm}

\begin{proof}
Let $\|\cdot\|$ be any norm on $\R^n$ and let $\lambda\in \R^{n*}$ be
a linear functional on $\R^n$ with $\|\lambda\|=1$.  Let $u\in \R^n$
be a vector with $\|u\|=1\lambda(u)$.  Let $S:=\{x\in
\R^n:\lambda(x)=0\}$ be the null space of $\R^n$.  Choose $n$~points
$a_0\cd a_n$ in $S$ that are affinely independent.  For each $M>0$
define 
$$
V_M:=\conv\{Ma_0\cd Ma_{n-1}\}.
$$
Any point of $V_M$ is uniquely of the form
$\sum_{k=0}^{n-1}x_kMa_k$ for some
$\sum_{k=0}^{n-1}x_ke_k\in\Delta_{n-1}$.  Define $F_M$ on $V_M$ by
$$
F_M\biggl(\,\sum_{k=0}^{n-1}x_kMa_k\biggr)
	=E\biggl(\,\sum_{k=0}^{n-1}x_ke_k\biggl).
$$
Finally set
$$
A_M:=\{x+yu:x\in V_M, F_M(x)\le y\le \kappa(n-1)+1\}.
$$
Since $E$ is lower semi-continuous $F_M$ is also lower
semi-continuous.  This implies $A_M$ is closed and bounded. (To see
$A_M$ is closed: $x_\ell+y_\ell u\in
A_M$ and $x_\ell+y_\ell u\to x+yu$ implies $x_\ell\to x$ and
$F_M(x)\le \liminf_{\ell\to\infty}F_M(x_\ell)\le
\lim_{\ell\to\infty}y_\ell=y$ and so $x+yu\in A_M$.)  It is also easy
to check $A_M$ is connected (and in fact contractible).  That $A_M$ is
an approximately convex sets follows from $E$ being an approximately convex
function.

Let $\e>0$ and define $\phi_M\cn\Delta_{n-1}\to V_M$ by 
$$
\phi_M(x)=\phi_M\biggl(\,\sum_{k=0}^{n-1}x_ke_k\biggr)=\sum_{k=0}^{n-1}x_kMa_k.
$$
Then $F_M\circ\phi_M=E$.  Fix a norm $\|\cdot\|_0$ on $\R^n$.  Then
there is a constant $C>0$ so that
$$
\|\phi_M(x)-\phi_M(y)\|\ge CM \|x-y\|_0\quad \text{for all}\quad
x,y\in \Delta_{n-1}.
$$
($C$ will depend on $\|\cdot\|_0$.) 
Since $E$ is lower semi-continuous $U:=\{x\in \Delta_{n-1}:
E(x)>E(\alpha)-\e\}$ is open in $\Delta_{n-1}$ and thus there is an
$R>0$ so that $B_R(\alpha)\cap \Delta_{n-1}\subset U$. Let
$a_k:=\phi_M(e_k)$.  Then as $E(e_k)=0$ we have $a_k\in A_M$ for $0\le
k\le n-1$. Let  $b:=\phi_M(\alpha)=\sum_{k=0}^{n-1}\alpha_ka_k$.
If $w\in A_M$ then $w=z+\beta u$ where $z\in
V_M$ and $F_M(z)\le \beta \le \kappa(n-1)+1$. If $\|z-b\| < MCR$, then
$F_M(z) > E(\alpha)-\e$ so that $\|z+\beta u-b\|\ge \lambda(\beta
u)=\beta \ge E(\alpha)-\e$. If $\|z-y\|\ge MCR$, then 
$$
\|z+\beta u-b\|\ge \|z-b\|-\beta\ge MCR -\kappa(n-1)-1.
$$
Now choose $M$ so that $MCR>2\kappa(n-1)+1$ so that $MCR
-\kappa(n-1)-1 \ge \kappa(n-1)\ge E(\alpha)-\e$.  
Then $\|z+\beta u-b\|\ge \kappa(n-1)-\e$ for all $z+\beta u\in A_M$
and so $\dist(b,A_M)\ge E(\alpha)-\e$.  This completes the proof.
\end{proof}

In the terminology of the last proof define a function $h_M\cn
\Delta_{n-1}\to [0,\infty)$ by $h_M:=\dist(\phi_M(x),A_M)$.  Then
$h_M$ is approximately convex and $h_M$ vanishes on the vertices of
$\Delta_{n-1}$.  Also $h_M$ is continuous and in fact Lipschitz
continuous.  The proof shows that for each fixed $\alpha\in
\Delta_{n-1}$ that $\lim_{M\to \infty}h_M(\alpha)=E(\alpha)$.
Replacing $n-1$ by $n$ we thus have:

\begin{prop}\label{E-cont-sup}
There is a sequence of Lipschitz continuous approximately convex functions
$\la h_\ell\ra_{\ell =0}^\infty$ on $\Delta_n$ vanishing on the
vertices of $\Delta_n$ such that $\lim_{\ell\to
\infty}h_\ell(x)=\sup_{\ell\ge 1}h_\ell(x)=E(x)$ for all $x\in
\Delta_n$.\qed
\end{prop}

\subsection{The sharp bounds in Euclidean Space}
\label{sec:euclidean}

Theorem~\ref{gen-upper} can be improved in Euclidean spaces.

\begin{thm} \label{thm:Euc-bd}
Let $\R^n$ have its usual inner product norm and let $A\subset \R^n$
be approximately convex.  Then any point $b\in \conv(A)$ has
$\dist(b,A)\le \kappa(n-1)$. (When combined with Theorem~\ref{II4} and
using the terminology of Theorem~\ref{midpoint} this implies 
$C_{\|\cdot\|}=\kappa(n-1)$ in Euclidean spaces of all dimensions.)
\end{thm}

We will denote the usual inner product on $\R^n$ by
$\la\cdot,\cdot\ra$.  Let $S^{n-1}$ be the unit sphere in $\R^n$ with
the Euclidean norm.  Set
\begin{equation}\label{Cn-def}
\simp(n):=\{A\subset S^{n-1}: \#(A)=n+1\ \text{and}\ 0\in\conv(A)^\circ\},
\end{equation}
so that $\simp(n)$ can be thought of as the set of simplexes inscribed
in the sphere that have the origin~$0$ in their interior.  An
$n$-dimensional simplex that has all its edge lengths equal is a {\bi
regular simplex\/}.  Recall that any two regular simplices with the
same edge lengths are congruent. We leave following calculations to
the reader.

\begin{prop}\label{reg:props}
Let $A\subset S^{n-1} $ be the set of vertices of a regular
$n$-dimensional simplex (so that $\#(A)=n+1$) inscribed in the
sphere.  Then $A\in \simp(n)$ and the  edge length of $A$ is given by
$$
\|a-b\|=\sqrt{\frac{2(n+1)}{n}}.
$$
($a,b\in A$ and $a\ne b$).  Moreover the distance of
the midpoint of the segment
between $a$ and $b$ to the origin is
$$
\biggl\|\frac{a+b}{2}\biggr\|= \sqrt{\frac{n-1}{2n}}.
$$
\vskip-15pt\qed
\end{prop}

Define $\maxside\cn \simp(n)\to [0,2]$ by
$$
\maxside(A)=\max_{a,b\in A}\|a-b\|.
$$
Then $\maxside(A)$ is the length of the longest edge of the simplex
with vertices $A$.  The following characterizes the regular simplexes
in terms of minimizing~$\maxside$ on $\simp(n)$.

\begin{thm}\label{reg-char}
Let $A\in \simp(n)$.  Then
$$
\maxside(A)\ge \sqrt{\frac{2(n+1)}{n}}
$$
with equality if and only if $A$ is the set of vertices of a regular
simplex.
\end{thm}

\begin{lemma}\label{lemma:smaller}
Let $A=\{x_0 \cd x_n\}\in \simp(n)$ and assume that
\begin{equation}\label{n-1<n}
\|x_0-x_{n-1}\|<\|x_0-x_n\|.
\end{equation}
Then there is a point $x_0^*\in S^{n-1}$ so that 
\begin{equation}\label{in-C(n)}
\{x_0^*,x_1,x_2\cd x_n\}\in \simp(n),
\end{equation}
\begin{equation}\label{smaller-n}
\|x_0^*-x_n\|<\|x_0-x_n\|,
\end{equation}
\begin{equation}\label{smaller-n-1}
\|x_0^*-x_{n-1}\|<\|x_0-x_n\|,
\end{equation}
\begin{equation}\label{same-i}
\|x_0^*-x_i\|=\|x_0-x_i\|,\quad 1\le i\le n-2.
\end{equation}
\end{lemma}

\begin{proof}
Since $0$ is in the interior of $\conv(A)$ any subset of $A$ of
size~$n$ will be linearly independent.  For $1\le i\le n$ define
$f_i\cn\R^n\to [0,\infty)$ and $\rho_i\cn S^n\to \R^n$ by
$$
f_i(x):=\|x-x_i\|,\quad \rho_i(x):=\|x-x_i\|.
$$
($\rho_i$ is the restriction of $f_i$ to $S^{n-1}$.)  Let $\nabla
f_i$ be the usual gradient of $f_i$ and $\nabla\rho_i$ the gradient of
$\rho_i$ as a function on $S^{n-1}$.  (That is $\nabla\rho_i$ is the
vector field tangent to $S^{n-1}$ so that for smooth curves $c(t)$ in
$S^{n-1}$ the equality $\frac{d}{dt}\rho_i(c(t))=\la
c'(t),\nabla\rho_i(c(t))\ra$ holds.)  Then a standard calculation
gives
$$
\nabla f_i(x)=\frac{x-x_i}{\|x-x_i\|}.
$$
As $\rho_i$ is the restriction of $f_i$ to $S^{n-1}$ the vector field
$\nabla \rho_i(x)$ is the orthogonal projection of $\nabla f_i(x)$ onto
the tangent space $T(S^{n-1})_x$ to $S^{n-1}$ at $x$. Therefore 
$$
\nabla\rho_i(x_0)=\frac{x_0-x_i}{\|x_0-x_i\|}-\Bigl\la
\frac{x_0-x_i}{\|x_0-x_i\|}, x_0 \Bigr\ra x_0.
$$
But then the $n-1$ vectors $\nabla\rho_1(x_0)\cd
\nabla\rho_{n-2}(x_0), \nabla\rho_n(x_0)$ are linearly independent as
any nontrivial linear relationship between them would lead to a
nontrivial linear relationship between $x_0\cd x_{n-2}, x_n$ which are
linearly independent.  The implicit function theorem implies that the
$n-1$ functions $\rho_1\cd \rho_{n-2},\rho_n$ are local coordinates on
$S^{n-1}$ near $x_0$ (that is the map $x\mapsto (\rho_1(x)\cd
\rho_{n-1}(x),\rho_n(x)$ is a diffeomorphism onto an open set in
$\R^{n-1}$ when restricted to a small enough open neighborhood of
$x_0$).  Let $\delta_i=\rho_i(x_0)=\|x_0-x_i\|$ and set
\begin{align*}
N:=&\{x\in S^{n-1}: \rho_i(x)=\delta_i, 1\le i\le n-2\}\\
=&\{x\in S^{n-1}: \|x-x_i\|=\delta_i, i\le i\le n-2 \}.
\end{align*}
As $\rho_1\cd \rho_{n-2},\rho_n$ are local coordinates near $x_0$ this
will be a smooth curve in $S^{n-1}$ near $x_0$ and, moreover, any point
$x_0^*\in N$ will satisfy all the conditions~(\ref{same-i}).  Choose a
parameterization $c\cn(-\e,\e)\to N$ of $N$ near $x_0$ with
$c(0)=x_0$.  As $\rho_1\cd \rho_{n-2},\rho_n$ is a local coordinate
system near $x_0$ and the first $n-2$ for these functions are constant
on $c(t)$ we have that $\frac{d}{dt}\rho_n(c(t))\big|_{t=0}
=\la\nabla\rho_n(x_0),c'(t)\ra\ne 0$.  Without loss of generality we
can assume that $\frac{d}{dt}\rho_i(c(t))\big|_{t=0}<0$ (otherwise
replace $c(t)$ by $c(-t)$).  Then for small $t>0$ we have
$\rho_n(c(t))<\rho_n(c(0))=\rho_n(x_0)$.  Also the
conditions~(\ref{in-C(n)} and~(\ref{smaller-n-1}) are open conditions
in $x_0^*$ and so for any $t$ sufficiently close to $0$ they will hold
for $x_0^*=c(t)$.  Therefore $x_0^*=c(t)$ for small positive $t$
satisfies the conclusion of the lemma.  This completes the proof.
\end{proof}

\begin{proof}[Proof of Theorem~\ref{reg-char}] We prove the theorem by
induction on $n$.  The base case of $n=1$ is trivial.  Let
$\ol{\simp}(n)$ be the closure of $\simp(n)$, that is
$$
\ol{\simp}(n)=\{A\subset S^{n-1}: \#(A)\le n+1, 0\in\conv(A)\}.
$$
Then the function $\maxside(A)=\max_{a,b\in A}\|a-b\|$ is continuous
on $\ol{\simp}(n)$ and $\ol{\simp}(n)$ is compact, so $\maxside$
obtains its minimum at some $A_0\in \ol{\simp}(n)$.  If this minimum
occurs at a boundary point of $\ol{\simp}(n)$ then $0\in \conv(A_0)$,
but $0\notin \conv(A_0)^\circ$.  Let $a,b\in A_0$ be the points of
$A_0$ so that $\|a-b\|=\maxside(A_0)$.  Then there exists a subset
$\{a,b\}\subseteq A_1\subseteq A_0$ so that $\#(A_1)=:m+1<n+1$ with
$A_1$ affinely independent and so that $0$ is in the relative interior
of $\conv(A_1)$.  Thus, with obvious notation, $A_1\in \simp(m)$ and
therefore by the induction hypothesis $\maxside(A_1)\ge
\sqrt{{2(m+1)}/{m}}$.  But for the regular simplex in $\simp(n)$ that
$\maxside$ has the value $\sqrt{{2(n+1)}/{n}}$, which is less than
$\sqrt{{2(m+1)}/{m}}$.  Therefore the minimum of $\maxside$ on
$\ol{\simp}(n)$ occurs in $\simp(n)$.

Again, let $A_0\in\simp(n)$ be where $\maxside$ obtains its minimum,
and let $c=\maxside(A_0)$.  If every edge of $A_0$ has length $c$ then
$A_0$ is a regular simplex and we are done.  The number of edges of
$A_0$ is $\binom{n+1}{2}$.  So assume that there are
$k<\binom{n+1}{2}$ edges that have length $c$.  Then there will be a
side $\{x_0,x_n\}$ of length $c$ that has a vertex in common with a
side $\{x_0,x_{n-1}\}$ that was a length less than $c$.  With this
notation let $A_0=\{x_0\cd x_n\}$.  Then by Lemma~\ref{lemma:smaller}
we can replace $x_0$ be some $x_0^*$ so that if $A_1:=\{x_0^*,x_1\cd
x_n\}$ then both the edges $\{x_0^*,x_n\}$ and $\{x_0^*,x_{n-1}\}$
have length $<c$ and all of the other $\binom{n+1}{2}-2$ edge lengths
stay the same.  Therefore $A_1$ has only $k-1$ edges of length $c$
(and if $k=1$ then all edges of $A_1$ have length less than $c$).  By
repeating this procedure $k$ times we end up with $A_k\in \simp(n)$ so
that $\maxside(A_k)<\maxside(A_0)$, contrary to the assumption that
$A_0$ was the minimizer.  Thus the minimizer must be regular.  This
completes the proof.
\end{proof}

If $x,y\in S^{n-1}$ then $\|x+y\|^2+\|x-y\|^2=4$.  Whence the distance
$\|\frac12(x+y)\|$ of the midpoint of the segment $\overline{xy}$ from
the origin is determined by its length.  Therefore
Theorem~\ref{reg-char} implies the following:

\begin{cor}\label{mid-char}
Let $\simp(n)$ be defined by~(\ref{Cn-def}) above and let
${\mathcal D}\cn \simp(n)\to [0,2]$ be given by
$$
{\mathcal D}(A)=\min_{a,b\in A}\biggl\|\frac{a+b}{2}\biggr\|.
$$
Then for all $A\in \simp(n)$ the inequality
$$
{\mathcal D}(A)\le \sqrt{\frac{n-1}{2n}}
$$
holds.  Equality holds if and only if $A$ is the set of vertices of a
regular simplex.\qed
\end{cor}

The following is what is needed in the proof of our main results.

\begin{prop}\label{mid<1}
Let $B_r(x_0)$ be a ball of radius $r$ in $\R^n$ with the Euclidean
norm and assume that there are $n+1$ points $\{a_0\cd a_n\}\subseteq
\partial B_r(x_0)$ such that $x_0$ is in the interior of the simplex
$\conv\{a_0\cd a_n\}$.  Assume that for each pair $\{a_i,a_j\}$  that
the distance of the midpoint $(a_i+a_j)/2$ to $\partial B_r(x_0)$ is
$\le 1 $.  Then 
$$
r\le \frac{\sqrt{2n}(\sqrt{2n}+\sqrt{n-1})}{n+1}\le \kappa(n-1) 
$$  
\end{prop}

\begin{proof}
By Corollary~\ref{mid-char} there exists a pair $\{a_i,a_j\}$ such
that $\|a_i+a_j\|/2\le r\sqrt{(n-1)/(2n)}$.  So
$$
1 \ge \dist\Big(\frac{a_i+a_j}2, \partial  B_r(x_0)\Big)\ge 
 r-\frac{\|a_i+a_j\|}{2}\ge r\(1-\sqrt{\frac{n-1}{2n}}\,\)
$$
Solving for $r$ gives $r\le \sqrt{2m}(\sqrt{2n}-\sqrt{n-1})/(n+1)$.
To see that $r=r(n)\le \kappa(n-1)$ first note $r(n)<2+\sqrt(2)<3.42$
for $n\ge 1$.  If $n\ge 4$ we then have $r(n)<3.5=\kappa(4-1)\le
\kappa(n-1)$.  This only leaves $r(2)=2=\kappa(2-1)$ and
$r(3)=\sqrt{3}(\sqrt{3}+1)/2<3=\kappa(3-1)$.
\end{proof}

\begin{proof}[Proof of Theorem~\ref{thm:Euc-bd}]  
By replacing $A$ by its closure we can assume that $A$ is closed.
Define $f\cn \R^n\to \R$ by $f(x):=\dist(x,A)$.  By Cara\-th\'eodory's
Theorem, it suffices to prove that if $\{a_0\cd a_n\}\subseteq A$ then
$f(x)\le \kappa(n-1)$ for all $x\in \conv(\{a_0\cd a_n\})$.  To
simplify notation set $\Delta:=\conv(\{a_0\cd a_n\})$ and let $x_0$ be
the point where $f\big|_{\Delta}$ achieves its maximum. Then we wish to
show $f(x_0)\le \kappa(n-1)$.  If $x_0$ is on the boundary (or if
$\{a_0\cd a_n\}$ is not affinely dependent) then $x_0$ is a convex
combination of $\le n$ points of $\{a_0\cd a_n\}$ and so $f(x_0)\le
\kappa(n-1)$ by Theorem~\ref{thm:convex-bd}.

This leaves the case where $x_0$ is in the interior of $\Delta$.  Then
$f(\cdot)=\dist(\cdot,A)$ has a local maximum at the interior point
$x_0$ of $\conv(A)$.  Let $R:=f(x_0)$.  Then, by
Proposition~\ref{trap0}, there are points $a_0\cd a_k\in A\cap
B_R(x_0)$ so that $\{a_0\cd a_k\}$ is an affinely independent set
and there are unit vectors $u_0\cd u_k$ so that the functional
$\lambda_i:=\la\cdot,u_i\ra$ norms $a_i-x_0$ and $0\in \conv\{u_0\cd
u_k\}$.  But if $\lambda_i$ norms $a_i-x_0$ then
$u_i=(a_i-x_i)/\|a_i-x_i\|$.  Therefore $0\in \conv\{u_0\cd u_k\}$
implies $x_0\in \conv\{a_0\cd a_k\}$.  Now Proposition~\ref{mid<1}
implies $f(x_0)=R\le \kappa(n-1)$.  This completes the proof.
\end{proof}

\begin{remark}\label{sharp-2D-Euc}
Let $A$ be the seven point subset of the Euclidean plane shown in
Figure~\ref{2D-ex}.  Then $A$ is approximately convex and satisfies
$d_H(\conv(A),A)=2=\kappa(1)$.  In higher dimensions we do not know if
there exist such examples of $A\subset \R^n$ with
 $d_H(\conv(A),A)=\kappa(n-1)$.
\begin{figure}[ht]
\centering
\setlength{\unitlength}{0.00033333in}
\begingroup\makeatletter\ifx\SetFigFont\undefined
\def\x#1#2#3#4#5#6#7\relax{\def\x{#1#2#3#4#5#6}}%
\expandafter\x\fmtname xxxxxx\relax \def\y{splain}%
\ifx\x\y   
\gdef\SetFigFont#1#2#3{%
  \ifnum #1<17\tiny\else \ifnum #1<20\small\else
  \ifnum #1<24\normalsize\else \ifnum #1<29\large\else
  \ifnum #1<34\Large\else \ifnum #1<41\LARGE\else
     \huge\fi\fi\fi\fi\fi\fi
  \csname #3\endcsname}%
\else
\gdef\SetFigFont#1#2#3{\begingroup
  \count@#1\relax \ifnum 25<\count@\count@25\fi
  \def\x{\endgroup\@setsize\SetFigFont{#2pt}}%
  \expandafter\x
    \csname \romannumeral\the\count@ pt\expandafter\endcsname
    \csname @\romannumeral\the\count@ pt\endcsname
  \csname #3\endcsname}%
\fi
\fi\endgroup
{\renewcommand{\dashlinestretch}{30}
\begin{picture}(13638,3339)(0,-10)
\path(12,312)(13212,312)
\path(7812,462)(7812,162)
\path(9012,462)(9012,162)
\path(10212,462)(10212,162)
\path(11412,462)(11412,162)
\path(12612,462)(12612,162)
\path(5487,462)(5487,162)
\path(5487,462)(5487,162)
\path(4212,462)(4212,162)
\path(3012,462)(3012,162)
\path(1812,462)(1812,162)
\path(612,462)(612,162)
\path(6462,2712)(6762,2712)
\path(6612,3312)(6612,12)
\path(6462,1437)(6762,1437)
\put(6537,1362){\makebox(0,0)[lb]{\smash{{{\SetFigFont{20}{24.0}{rm}.}}}}}
\put(8637,2637){\makebox(0,0)[lb]{\smash{{{\SetFigFont{20}{24.0}{rm}.}}}}}
\put(4587,2637){\makebox(0,0)[lb]{\smash{{{\SetFigFont{20}{24.0}{rm}.}}}}}
\put(10737,1437){\makebox(0,0)[lb]{\smash{{{\SetFigFont{20}{24.0}{rm}.}}}}}
\put(2412,1437){\makebox(0,0)[lb]{\smash{{{\SetFigFont{20}{24.0}{rm}.}}}}}
\put(12537,312){\makebox(0,0)[lb]{\smash{{{\SetFigFont{20}{24.0}{rm}.}}}}}
\put(537,312){\makebox(0,0)[lb]{\smash{{{\SetFigFont{20}{24.0}{rm}.}}}}}
\put(10962,1512){\makebox(0,0)[lb]{\smash{{{\SetFigFont{7}{8.4}{rm}$(2\sqrt3,1)$}}}}}
\put(2637,1512){\makebox(0,0)[lb]{\smash{{{\SetFigFont{7}{8.4}{rm}$(-2\sqrt3,1)$}}}}}
\put(8862,2712){\makebox(0,0)[lb]{\smash{{{\SetFigFont{7}{8.4}{rm}$(\sqrt3 ,2)$}}}}}
\put(4812,2712){\makebox(0,0)[lb]{\smash{{{\SetFigFont{7}{8.4}{rm}$(-\sqrt3 ,2)$}}}}}
\put(6912,1437){\makebox(0,0)[lb]{\smash{{{\SetFigFont{7}{8.4}{rm}$(0,1)$}}}}}
\put(162,612){\makebox(0,0)[lb]{\smash{{{\SetFigFont{7}{8.4}{rm}$(-3\sqrt3,0)$}}}}}
\put(12162,612){\makebox(0,0)[lb]{\smash{{{\SetFigFont{7}{8.4}{rm}$(3\sqrt3,0)$}}}}}
\end{picture}
}
\caption[]{A two-dimensional Euclidean example.}
\label{2D-ex}
\end{figure}
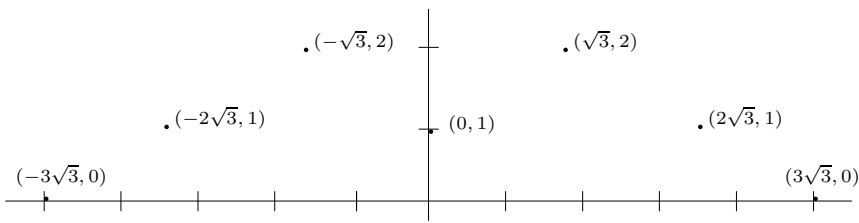
\end{remark}

\subsection{The sharp two dimensional bounds}
\label{sec:2D}
We now give the sharp estimate for the size of a convex hull in all
two dimensional normed spaces.
\begin{thm}\label{2D-bds}
Suppose $\|\cdot\|$ is a norm on $\R^2$ and that $A \subseteq X$ has
is approximately convex in this norm. Then any point $b\in \conv(A)$ has
$\dist(b,A)\le 2$.  (By Theorem~\ref{II4} given $\e>0$, there exists
an approximately convex $A_\e \subseteq \R^2$ and a $b\in \conv(A)$ so that
$\dist(b,A)\ge 2-\e$ and thus thus in the notation of
Theorem~\ref{midpoint} $C_{\|\cdot\|}=2$ for all two dimensional norms.)
\end{thm}

\begin{lemma}\label{steve1} 
Let $V=\{a,b,c,-a,-b,-c\}$ be the
vertices of a symmetric convex hexagon. Then
$$\{a+b,b+c,c+a\}\cap \conv(V) \ne \nothing.
$$ 
\end{lemma}

\begin{proof} By applying a linear transformation we may assume 
$a=(-1,1)$ and $b=(-1,-1)$.  Without loss
of generality we also assume
$c=(x_0,y_0)$, where $-1 \le y_0 \le 0$ and $x_0\ge 1$.
If $y_0 > 2 - x_0$, then $a+b=(-2,0) \in \conv(V)$, and
we are done. So we may assume that 
$ c \in \conv(\{(1,0),(1,-1),(2,0),(3,-1)\})$.
\begin{figure}[ht]
\centering
\setlength{\unitlength}{0.00037500in}
\begingroup\makeatletter\ifx\SetFigFont\undefined
\def\x#1#2#3#4#5#6#7\relax{\def\x{#1#2#3#4#5#6}}%
\expandafter\x\fmtname xxxxxx\relax \def\y{splain}%
\ifx\x\y   
\gdef\SetFigFont#1#2#3{%
  \ifnum #1<17\tiny\else \ifnum #1<20\small\else
  \ifnum #1<24\normalsize\else \ifnum #1<29\large\else
  \ifnum #1<34\Large\else \ifnum #1<41\LARGE\else
     \huge\fi\fi\fi\fi\fi\fi
  \csname #3\endcsname}%
\else
\gdef\SetFigFont#1#2#3{\begingroup
  \count@#1\relax \ifnum 25<\count@\count@25\fi
  \def\x{\endgroup\@setsize\SetFigFont{#2pt}}%
  \expandafter\x
    \csname \romannumeral\the\count@ pt\expandafter\endcsname
    \csname @\romannumeral\the\count@ pt\endcsname
  \csname #3\endcsname}%
\fi
\fi\endgroup
{\renewcommand{\dashlinestretch}{30}
\begin{picture}(6815,4839)(0,-10)
\texture{55888888 88555555 5522a222 a2555555 55888888 88555555 552a2a2a 2a555555 
	55888888 88555555 55a222a2 22555555 55888888 88555555 552a2a2a 2a555555 
	55888888 88555555 5522a222 a2555555 55888888 88555555 552a2a2a 2a555555 
	55888888 88555555 55a222a2 22555555 55888888 88555555 552a2a2a 2a555555 }
\put(1257,3597){\shade\ellipse{94}{94}}
\put(1257,3597){\ellipse{94}{94}}
\put(1257,1197){\shade\ellipse{94}{94}}
\put(1257,1197){\ellipse{94}{94}}
\put(3657,1197){\shade\ellipse{94}{94}}
\put(3657,1197){\ellipse{94}{94}}
\put(3657,3597){\shade\ellipse{94}{94}}
\put(3657,3597){\ellipse{94}{94}}
\put(357,2997){\shade\ellipse{94}{94}}
\put(357,2997){\ellipse{94}{94}}
\put(4557,1797){\shade\ellipse{94}{94}}
\put(4557,1797){\ellipse{94}{94}}
\put(4859,2412){\shade\ellipse{94}{94}}
\put(4859,2412){\ellipse{94}{94}}
\put(6059,1212){\shade\ellipse{94}{94}}
\put(6059,1212){\ellipse{94}{94}}
\put(3659,2412){\shade\ellipse{94}{94}}
\put(3659,2412){\ellipse{94}{94}}
\put(3357,2997){\shade\ellipse{94}{94}}
\put(3357,2997){\ellipse{94}{94}}
\path(12,2412)(6012,2412)
\path(2412,4812)(2412,12)
\path(1212,1212)(3612,1212)(4512,1812)
	(3612,3612)(1212,3612)(312,3012)(1212,1212)
\path(3612,1212)(6012,1212)(4812,2412)
	(3612,2412)(3612,1212)
\dashline{45.000}(5112,2112)(5112,1212)
\dashline{45.000}(5412,1812)(5412,1212)
\dashline{45.000}(5712,1512)(5712,1212)
\dashline{45.000}(4812,2412)(4812,1212)
\dashline{45.000}(4512,2397)(4512,1197)
\dashline{45.000}(4212,2412)(4212,1212)
\dashline{45.000}(3912,2412)(3912,1212)
\dashline{45.000}(3612,2112)(5112,2112)
\dashline{45.000}(3612,1812)(5412,1812)
\dashline{45.000}(3612,1512)(5712,1512)
\path(3537,3012)(5112,3312)
\path(3649.267,3063.924)(3537.000,3012.000)(3660.494,3004.983)
\put(3237,2562){\makebox(0,0)[lb]{\smash{{{\SetFigFont{8}{9.6}{rm}$(1,0)$}}}}}
\put(4662,2562){\makebox(0,0)[lb]{\smash{{{\SetFigFont{8}{9.6}{rm}$(2,0)$}}}}}
\put(537,2862){\makebox(0,0)[lb]{\smash{{{\SetFigFont{8}{9.6}{rm}$-c$}}}}}
\put(4587,1587){\makebox(0,0)[lb]{\smash{{{\SetFigFont{8}{9.6}{rm}$c$}}}}}
\put(612,912){\makebox(0,0)[lb]{\smash{{{\SetFigFont{8}{9.6}{rm}$b=(-1,-1)$}}}}}
\put(687,3687){\makebox(0,0)[lb]{\smash{{{\SetFigFont{8}{9.6}{rm}$a=(-1,1)$}}}}}
\put(3387,3687){\makebox(0,0)[lb]{\smash{{{\SetFigFont{8}{9.6}{rm}$-b=(1,1)$}}}}}
\put(3312,912){\makebox(0,0)[lb]{\smash{{{\SetFigFont{8}{9.6}{rm}$-a=(1,-1)$}}}}}
\put(5862,912){\makebox(0,0)[lb]{\smash{{{\SetFigFont{8}{9.6}{rm}$(3,-1)$}}}}}
\put(5187,3312){\makebox(0,0)[lb]{\smash{{{\SetFigFont{8}{9.6}{rm}$a+c$}}}}}
\end{picture}
}
\caption[]{}
\label{hex}
\end{figure}
($\conv(\{(1,0),(1,-1),(2,0),(3,-1)\})$ is shaped region in
Figure~\ref{hex}.)  
This forces the quadrilateral
$\conv(\{0,a,c,-b\})$ to contain the parallelogram
$\conv(\{0,a,c,a+c\})$, and so
$$a+c \in \conv(\{0,a,c,a+c\})\subseteq
 \conv(\{0,a,c,-b\}) \subseteq
\conv(V).$$ 
\vskip-15pt
\end{proof}

For the rest of this section we will call a norm on a finite
dimensional space $\|\cdot\|$ {\bi smooth\/} if it is a $C^\infty$
function away from the origin and the unit ball is strictly convex.  A
finite dimensional space is smooth iff its norm is smooth.  This
implies that norming linear functionals are unique.

\begin{lemma} \label{steve2}
Let $X$ be a smooth two-dimensional normed space.
Suppose that $K \subseteq S_1(0)$ is a closed set and that
$0 \notin \conv(K)$. Then $f(x)=\dist(x,K)$ does not attain a local
maximum at $x=0$. 
\end{lemma}

\begin{proof} 
As $(X,\|\cdot\|)$ is smooth for each $u\in S_1(0)$ there is a unique
norm linear functional $\lambda_u$ that norms $u$, the map $u\mapsto
\lambda_u$ is a homeomorphism of $S_1(0)$ onto the unit sphere
$S^*_1(0)$ in the dual space $(X^*,\|\cdot\|^*)$, and
$\lambda_{-u}=-\lambda_u$.  If $u\in S_1(0)$ then
$S_1(0)\setminus\{u,-u\}$ has exactly two connected components.  A
closed subset $K\subseteq S_1(0)$ satisfies $0\notin \conv(K)$ if and
only if there is a $u\in S_1(0)$ so that $K$ is contained in one of
the connected components of $S_1(0)\setminus \{u,-u\}$ (for this is
equivalent to being able to separate $K$ from the origin by a linear
functional).  But the properties of the map $u\mapsto \lambda_u$ imply
$K$ is contained in a connected component of $S_1(0)\setminus
\{-u,u\}$ if and only if $N^*(K):=\{\lambda_u: u\in K\}$ is contained
in a connected component of
$S_1^*(0)\setminus\{\lambda_u,-\lambda_u\}$.  Therefore $0\notin
\conv(K)$ if and only if $0\notin \conv(N^*(K))$.  But by
Proposition~\ref{trap0} $0\notin \conv(N^*(K))$ implies that $f$ does
not have a local maximum at $0$.
\end{proof}

Let $\e>0$.
A set $A\subseteq X$ will be said to be {\bi $\e$-separated\/}
if $\|a-b\|\ge \e$ whenever $a,b$ are distinct
elements of $A$.
\begin{lemma}\label{steve3}
Suppose that $X$ is a smooth two-dimensional normed space and that $A
\subseteq X$ is $\e$-separated and approximately convex. Then
$d_H(A,\conv(A))\le 2$.
\end{lemma}
\begin{proof} Let $f(x) = \dist(x,A)$ ($x \in X$). By Carath\'eodory's
Theorem, it suffices to prove that if $\{d,e,f\} \subseteq A$, then
$f(x) \le 2$ for all $x \in \Delta$, where $\Delta = \conv(
\{d,e,f\})$. By continuity of $f$, there exists $x_0 \in \Delta$ at
which $f$ attains its maximum. By translation we may assume without
loss of generality that $x_0=0$. If $0 \in \partial( \Delta)$ then $0$
is on a segment between two elements of $A$ and so by restriction $f$
to this segment see by Theorem~\ref{thm:convex-bd} $f(0) \le 2$. So we
may assume that $0$ lies in the interior of $\Delta$.  Let $R=f(0)$
and let $K=A \cap B_R(0)$. If $0 \notin \conv(K)$, then by
Lemma~\ref{steve2} $g(x)= \dist(x,K)$ does not attain a local maximum
at $x=0$.  But since $A$ is $\e$-separated an easy compactness
argument yields $\dist(0,A\setminus K) >R$, and so $f(x)=g(x)$ for all
$x$ sufficiently close to $x=0$. Thus, $f(x)$ does not attain a local
maximum at $x=0$, which contradicts the fact that $0$ lies in the
interior of $\Delta$.

So we may assume that $0\in \conv(K)$. By Carath\'eodory's  Theorem
there exists $\{a,b,c\}\subseteq K$ with $0 \in \conv(\{a,b,c\})$
Once again, we may assume that $0$ lies in the interior
of $\conv(\{a,b,c\})$. Now $B_R(0)$ contains the convex
hexagon with vertices $V=\{a,b,c,-a,-b,-c\}$. By Lemma~\ref{steve2},
$$\{a+b,b+c,c+a\}\cap \conv(V) \ne \emptyset.$$ Thus
$$\min\{\|a+b\|,\|b+c\|,\|c+a\|\}\le R.$$
We may assume without loss of generality that $\|a+b\|\le R$.
Since $A$ is approximately convex there exists $x \in A$
with $\|x-(1/2)(a+b)\| \le 1$. Thus
$$
R = \dist(0,A) \le \|x\| \le 1 + \frac{1}{2}\|a+b\| \le 1+
\frac{R}{2},
$$
and so $R \le 2$ as required. 
\end{proof}

\begin{proof}[Proof of Theorem~\ref{2D-bds}] 
Assume $A \subseteq X$ is approximately convex Let
$\e>0$.  There exists an equivalent smooth norm $\|\cdot\|'$ on $X$
such that
$$\|x\|' \le \|x\| \le (1+\e)\|x\|' \qquad (x \in  X).$$
Let $B \subseteq A$ be a maximal $\e$-separated subset of $A$.
Then $d_H'(A,B) \le d_H(A,B) \le \e$ (here $d_H'(\cdot,\cdot)$
denotes Hausdorff distance with respect to $\|\cdot\|'$).
Thus,
\begin{align*}
d_H'\left(B,\frac{B+B}{2}\right) 
& \le  	 d_H'(B,A)+ d_H'\left(A,\frac{A+A}{2}\right) 
+ d_H'\left(\frac{A+A}{2},\frac{B+B}{2}\right)\\
& \le  \e + 1 + \e =1 + 2 \e
\end{align*}  Lemma~\ref{steve3}  applied to $\|\cdot\|'$ and $B$
yields $d_H'(B,\conv(B)) \le 2(1+2\e)$. Thus,
\begin{align*}
d_H(A,\conv(A)) & \le  (1+\e)d_H(A,\conv(A)) \\
& \le  (1+\e)((d_H'(B,\conv(B))+ 2d_H'(A,B))\\
& \le  (1+\e)(2(1+2\e)+2\e)
\end{align*} Since $\e>0$ is arbitrary, we obtain
$d_H(A,\conv(A)) \le 2$ as desired.
\end{proof}

{\small {\bf Acknowledgments:} We would like to thank Maria Girardi
for realizing a blackboard on the North end of the third floor of the
Mathematics building would lead to the type of results given here.
During all stages of this work we profited from conversations with Anton
Schep.
}


\providecommand{\bysame}{\leavevmode\hbox to3em{\hrulefill}\thinspace}

\end{document}